\documentclass[a4paper,10pt,reqno]{elsarticle-mr}

\usepackage[%
   colorlinks=true,
   urlcolor=rltblue,       
   filecolor=rltgreen,     
   citecolor=rltgreen,     
   linkcolor=rltred,       
   bookmarks=true,
   plainpages=false,
   ]{hyperref}

\usepackage[operator,rose,frak,novargreek]{paper_diening}  
\usepackage[active]{srcltx}
\usepackage{amsmath,color}
\usepackage{amssymb,amsthm,amsxtra,amscd}

\usepackage{paralist}

\usepackage{mathrsfs} 

\definecolor{rltred}{rgb}{0.75,0,0}
\definecolor{rltgreen}{rgb}{0,0.5,0}
\definecolor{rltblue}{rgb}{0,0,0.75}

\newtheorem{Def}[equation]{Definition}
\newtheorem{Sa}[equation]{Theorem}
\newtheorem{Lem}[equation]{Lemma}

\newtheorem{Bem}[equation]{Remark}
\newtheorem{Kor}[equation]{Corollary}
\newtheorem{Vss}[equation]{Assumption}

\newenvironment{Bew}{\begin{proof}[Proof]}{\end{proof}}

\newcommand{\Nat}{\mathbb{N}}

\newcommand{\R}{\mathbb{R}}
\newcommand{\E}{\mathbf{E}}

\newcommand{\w}{\boldsymbol\omega}

\newcommand{\x}{\boldsymbol\xi}

\def\Xint#1{\mathchoice
   {\XXint\displaystyle\textstyle{#1}}%
   {\XXint\textstyle\scriptstyle{#1}}%
   {\XXint\scriptstyle\scriptscriptstyle{#1}}%
   {\XXint\scriptscriptstyle\scriptscriptstyle{#1}}%
   \!\int}
\def\XXint#1#2#3{{\setbox0=\hbox{$#1{#2#3}{\int}$}
     \vcenter{\hbox{$#2#3$}}\kern-.5\wd0}}
\def\dint{\Xint-}

\DeclareMathOperator{\Div}{\operatorname{div}}

\newcommand{\anti}{{\ensuremath{\mathrm{skew}}}}
\newcommand{\bell}{\boldsymbol{\ell}}

\newcommand{\lmk}{\lambda_{m,k}}
\newcommand{\lMk}{\lambda_{M,k}}
\newcommand{\omk}{\mathcal{O}_{m,k}}

\newcommand{\amk}{\alpha_{m,k}}
\newcommand{\aMk}{\alpha_{M,k}}
\newcommand\fett[1]{\boldsymbol{#1}}
\newcommand\vek[1]{\mathbf{#1}}
\DeclareMathOperator{\trunc}{{\mathcal T}^\alpha_E}
\DeclareMathOperator{\trunckm}{{\mathcal T}_{M,k}}

\newcommand{\omt}{\Omega_T}
\newcommand{\emk}{E_{M,k}}
\newcommand\grs{\vek{S}}
\newcommand\grw{\vek{W}}
\newcommand\grd{\vek{D}}
\newcommand\grr{\vek{R}}
\newcommand\gre{\vek{E}}
\newcommand\grn{\vek{N}}

\newcommand\tphi{\fett{\varphi}}
\newcommand\tpsi{\fett{\psi}}

\newcommand{\om}{\boldsymbol{\omega}}
\newcommand\vvv{\vek{v}}
\newcommand\f{\vek{f}}
\renewcommand\l{\fett{\ell}}
\newcommand{\sdre}{\grs\big(\vek{D}\vvv,\grr(\vvv,\om),\gre\big)}
\newcommand{\nome}{\grn(\nabla\om,\gre)}
\newcommand{\dv}{\vek{D}\vvv}

\newcommand\uuum{\vek{u}_M}
\newcommand{\vqom}{V_q(\Omega)}
\newcommand{\lzweiom}{L^2(\Omega)}
\newcommand{\hhom}{H(\Omega)}
\newcommand{\wqom}{W_0^{1,q}(\Omega)}

\newcommand{\levy}{\fett{\varepsilon}}

\numberwithin{equation}{section}

\begin{document}
\begin{frontmatter}

\title{Existence of weak solutions for unsteady motions of micro-polar electrorheological fluids}

\author[eb]{E.~B{\"a}umle}
\ead{Erik.Baeumle@gmx.de}

\author[mr]{M.~R\r u\v zi\v cka \corref{cor1}}
\ead{rose@mathematik.uni-freiburg.de}

\cortext[cor1]{Corresponding author}

\address[eb]{Bl{\"u}tenweg 12, D-77656 Offenburg, GERMANY} 

\address[mr]{Institute of Applied Mathematics,
  Albert-Ludwigs-University Freiburg, Eckerstr.~1, D-79104 Freiburg,
  GERMANY.}

\begin{abstract}
  In this paper we study the existence of weak solutions to an
  unsteady system describing the motion of micro-polar
  electrorheological fluids. The constitutive relations for the stress
  tensors belong to the class of generalized Newtonian fluids. Using
  the Lipschitz truncation and the solenoidal Lipschitz truncation we
  establish the existence of global solutions for shear exponents
  $p>6/5$ in three-dimensional domains.
\end{abstract}


\end{frontmatter}

\section{Introduction}\label{introduction}
In this paper we investigate the existence of solutions of the
system\footnote{We denote by $\levy$ the isotropic third order tensor
  and by $\levy : \grs$ the vector having the components
  $\varepsilon_{ijk}S_{jk}$, $i=1,\ldots d$, where the summation
  convention is used. }
\begin{align}
  \begin{aligned}\label{system}
    \partial_t \bfv + \operatorname{div}(\bfv\otimes\bfv) -
    \operatorname{div}\grs + \nabla \pi&=\bff &&\text{in} \ \Omega_T\,,
    \\
    \operatorname{div}\bfv &=0 &&\text{in} \ \Omega_T\,,
    \\
    \partial_t \om + \operatorname{div}(\om\otimes\bfv) -
    \operatorname{div}\grn &=\l -\levy:\grs &&\text{in} \ \Omega_T\,,
  \end{aligned}
\end{align}
completed with homogeneous Dirichlet boundary conditions
\begin{align}\label{bc}
  \mathbf{v}&=\bfzero\,, \quad \boldsymbol\bfomega =\bfzero \qquad
  \text{on}\ I \times \partial\Omega \,,
\end{align}
and initial conditions
\begin{align}\label{ic}
  \mathbf{v}(0)&=\bfv_0\,, \quad \boldsymbol\bfomega(0) =\bfomega _0 \qquad
  \text{in}\ \Omega \,.
\end{align}
Here $\Omega\subset \setR^3$ is a bounded domain and $I=(0,T)$ with $T
\in (0,\infty)$ a given finite time intervall. The three
equations in \eqref{system} are the balance of momentum, mass and angular
momentum for an incompressible, micro-polar electrorheological
fluid.  In these equations $\bv$ denotes the velocity, $\w$ the
micro-rotation, $\pi$ the pressure, $\bS$ the mechanical extra
stress tensor, $\bN$ the couple stress tensor, $\bell$ the
electromagnetic couple force, ${\ff=\tilde \ff + \chi^E\divo (\bE
  \otimes \bE)}$ the body force, where $\tilde \ff$ is the mechanical
body force, $\chi^E$ the dielectric susceptibility and $\bE$ the
electric field. The electric field $\bE$ solves quasi-static
Maxwell's equations
\begin{align}
  \begin{aligned}\label{maxwell}
    \Div \mathbf{E}&=0 &&\text{in}\ \Omega_T\,,
    \\
    \curl \mathbf{E}&=\bfzero &&\text{in}\ \Omega_T\,,
    \\
    \mathbf{E}\cdot \mathbf{n}&=\mathbf{E}_0\cdot \mathbf{n}
    &\quad&\text{on}\ I \times \partial\Omega,
  \end{aligned}
\end{align}
where $\mathbf{n}$ is the outer normal vector of the boundary
$\partial \Omega$ and $\mathbf{E}_0$ is a given time-dependent
electric field. The model \eqref{system}--\eqref{maxwell} is derived
in \cite{win-r}. It contains a more realistic description of the
dependence of the electrorheological effect on the direction of the
electric field compared to the previous model in \cite{RR2},
\cite{rubo}. Nevertheless, we concentrate in this paper on the
investigation of the mechanical properties of electrorheological
fluids governed by \eqref{system}. This is possible due to the fact
that Maxwell's equations \eqref{maxwell} are separated from the
balance laws \eqref{system} and that there exists a well developed
existence theory for Maxwell's equations. Thus, we will assume
throughout the paper that an electric field $\bE$ with appropriate
properties is given (cf.~Assumption~\ref{VssE}).

A representative example for a constitutive relation for the stress
tensors in \eqref{system} reads (cf.~\cite{win-r}, \cite{rubo})
\begin{align}
  \hspace*{-1mm}  
  \begin{aligned}\label{eq:SN-ex}
    \mathbf{S}&=(\alpha_{31}+\alpha_{33}\vert\mathbf{E}\vert^2)
    (1+\vert\mathbf{D}\vert)^{p-2}\mathbf{D}+ \alpha_{51}
    (1+\vert\mathbf{D}\vert)^{p-2}\big (\mathbf{D} \bE \otimes \bE +
    \bE \otimes \bD\bE \big)\hspace*{-5mm}
    \\
    &\quad +
    \alpha_{71}\vert\mathbf{E}\vert^2(1+\vert\mathbf{R}\vert)^{p-2}
    \mathbf{R} + \alpha_{91} (1+\vert\mathbf{R}\vert)^{p-2}\big
    (\mathbf{R} \bE \otimes \bE + \bE \otimes \bR\bE \big)\,,
    \\
    \mathbf{N}&=(\beta_{31}+\beta_{33}\vert\mathbf{E}\vert^2)
    (1+\vert\nabla\boldsymbol\omega\vert)^{p-2}\nabla\boldsymbol\omega
    \\
    &\quad + \beta_{51}(1+\vert\nabla \w\vert)^{p-2}\big ((\nabla \w
    )\bE \otimes \bE + \bE \otimes (\nabla \w)\bE \big)\,,
  \end{aligned}
\end{align}
with constants $\alpha_{31},\alpha_{33},\alpha_{71},\beta_{33}>
0$ and $\beta_{31}\ge 0$. The constants $\alpha_{51}, \alpha_{91}, \beta_{51}$
have to satisfy certain restrictions (cf.~\cite{win-r}, \cite{rubo}),
which ensure the validity of the second law of thermodynamics. In
\eqref{eq:SN-ex} we used the notation\footnote{Here $\bfvarepsilon\cdot\om $
  denotes the tensor with components $\vep_{ijk}\omega_k$,
  $i,j=1,\ldots,d$.} $\bD = (\nabla \bv)^\sym$, ${\bR= \bW \bv
  +\bfvarepsilon \cdot\om }$, with $\bW\bv =(\nabla \bv)^\anti$. In the
present paper we refrain from considering concrete constitutive
relations for the stress tensors, but we make general assumptions
covering prototypical situations (cf.~Assumption \ref{VssS} and
Assumption \ref{VssN}).

Micro-polar fluids have been introduced by Eringen in the sixties
(cf.~\cite{eringen-book} for an exhaustive treatment).
Electrorheological fluids can be modelled in various ways, see
e.g.~\cite{atbul}, \cite{RW}, \cite{winr1}, \cite{RR2},
\cite{win-r}. While there exists many investigations of micro-polar as
well as of electrorheological fluids (cf.~\cite{Lukaszewicz},
\cite{rubo}), there exists to our knowledge no investigations of
micro-polar electrorheological fluids except \cite{erw}, which is
based on the PhD thesis \cite{frank-phd} and the diploma thesis
\cite{weber-dipl}; and the diploma thesis \cite{erik-dipl}. The
present paper is based on the latter thesis.

\smallskip In the next section we introduce the notation, the
functional setting, give assumptions for the stress tensors and
collect some auxiliary results. In particular, the properties of the
solenoidal unsteady Lipschitz truncation are stated and a
generalization of the unsteady Lipschitz truncation is discussed. In
Section~\ref{seceasy} we present the analysis of our problem in the
context of pseudomonotone operator theory, which applies for shear
exponents $p \ge 11/5$. With the same tools we construct approximate
solutions in the more interesting case $p <11/5$ in
Section~\ref{secapprox}. Using the different Lipschitz truncations
we prove our main result in Section~\ref{secmain}.

\section{Preliminaries}\label{veroeffentlichung1}

\subsection{Notation and function spaces}\label{subsecnotation}
We denote by $c$ generic constants, which may change from line to
line. Scalar-valued functions will be written in normal font,
e.g.,~$f,\zeta$ while vector-valued functions will be denoted by
boldfaced letters, e.g.,~$\vek{u},{\bfvarphi}$. Capital boldface
letters will be used for tensor-valued functions\footnote{The only
  exception will be the electric field which is denoted as usual by
  $\gre$.}, e.g.,~$\grs$. The standard scalar product for vectors is
denoted by $\vek{v}\cdot \om$, while the standard scalar product for
tensors is denoted by $\vek{A}:\vek{B}$.  We use the usual Lebesgue
and Sobolev spaces $L^p(\Omega),\ W^{k,p}(\Omega)$,
$1\leq p\leq\infty$, $k \in \Nat$, where $\Omega \subset \setR ^3$ is
bounded domain with Lipschitz-boundary. For given $T \in (0,\infty)$
we use the notation $\Omega_T:= (0,T)\times \Omega$. By
$W_0^{1,p}(\Omega)$ we denote the completion of $C_0^\infty(\Omega)$
in $W^{1,p}(\Omega)$ which will be equipped with the gradient norm
$\norm{\nabla\, \cdot\, }_{L^p}$.  In the notation of function spaces
we do not distinguish between scalar, vector-valued and tensor-valued
spaces.  For $\vek{v} \in W^{1,p}(\Omega)$ we denote by $\grd\vek{v}$
the symmetric and by $\grw\vek{v}$ the skew-symmetric part of the
gradient,
i.e.~$ \grd\vek{v}= \frac{1}{2}( \nabla \vek{v} +\nabla \vek{v}^\top)$
and
$ \grw\vek{v}= \frac{1}{2}( \nabla \vek{v} -\nabla \vek{v}^\top )$.
Further we define for vectors $\vek{v}, \vek{\om}$ the
tensor 
$\grr(\vek{v},\om):= \grw\vek{v}+\levy\cdot \om$.  These definitions
imply the equality
\begin{align}\label{rulelevy}\begin{split}
    \grs:\nabla\vvv+(\levy:\grs)\cdot
    \om
    =\grs:\big(\grd\vvv + \grr(\vvv,\om) \big).
\end{split}\end{align}
Moreover,  for any $\om \in\R^d$, $\grs \in \R^{d\times d}$ there holds
\begin{align}\begin{split}\label{levyabsch}
\abs{\levy:\grs}\leq c\,\abs{\grs}, \quad
\abs{\levy\cdot\om}\leq c\,\abs{\om}.
\end{split}
\end{align}
Additionally we need some completions of
$\mathcal{V} (\Omega) := \lbrace \vek{u} \in C_0^\infty(\Omega) \fdg
\operatorname{div} \vek{u} =0\rbrace$.
We define $H(\Omega) := \overline{\mathcal{V} (\Omega)}^{L^2}, $
$V_p(\Omega) := \overline{\mathcal{V} (\Omega)}^{W^{1,p}},$
$V^3(\Omega) := \overline{\mathcal{V} (\Omega)}^{W^{3,2}}.$ While we
will use the usual $L^2$-norm on $H(\Omega)$ and the usual
$W^{3,2}$-norm on $V^3(\Omega)$, we will equip the space $V_p(\Omega)$
for $1<p<\infty$ with the norm
$\norm{\,\cdot\,}_{V_p}:= \norm{\grd\, \cdot\,}_{L^p}$, which defines
an equivalent norm due to Korn's inequality (cf.~\cite{boy-fab}). The
duality pairing between a Banach space $V$ and its dual $V^\ast$ will
be denoted by $\langle \cdot,\cdot\rangle_V$. 
We use the usual notation for Bochner spaces (cf.~\cite{GGZ},
\cite{zei-IIA}) and denote by $\frac{du}{dt}$ the generalized
derivative, i.e.~let $V,W$ be Banach spaces with a dense embedding
$V\hookrightarrow W$ and assume that for $u \in L^p(0,T;V)$ there
exists $w\in L^q(0,T;W)$ such that
$\int_0^T \varphi'(t)u(t)\,dt = -\int_0^T \varphi(t)w(t)\,dt $ holds
for any $\varphi \in C_0^\infty((0,T))$, then we set
$\frac{du}{dt}:=w$.  We introduce the Bochner-Sobolev space
\[
  W^{1,p,q}(0,T;V,W):= \Big\lbrace u \in L^p(0,T;V) \fdg
    \frac{du}{dt} \in L^q(0,T;W) \Big\rbrace.
\]
For details concerning these spaces we refer to \cite{boy-fab}. Let $(V,H,V^\ast)$ be a Gelfand-triple, 
i.e.~$V$ is a Banach space and $H$ is a Hilbert space, such that $V$ embeds densely into $H$. Then we define 
the Bochner-Sobolev space
\[
  W^1_p(0,T;V,H):= \Big\lbrace u \in L^p(0,T;V) \fdg \frac{du}{dt} \in
  L^{p'}(0,T;V^\ast) \Big\rbrace.
\]
It is well known that (cf.~\cite{boy-fab}, \cite{zei-IIA}) we have the continuous embeddings
\begin{align}\begin{split}\label{einbettung1}
W^{1,p,q}(0,T;V,W)&\hookrightarrow C(0,T;W), \\
W^1_p(0,T;V,H)&\hookrightarrow C(0,T;H),
\end{split}\end{align} 
and that for $u,v \in W^1_p(0,T;V,H)$ there holds the integration by parts formula, i.e.~for any 
$0\leq s,t \leq T$ there holds (cf.~\cite{zei-IIA})
\begin{align}\begin{split}\label{partielleInteg}
(u(t),v(t))_H - (u(s),v(s))_H = \int_s^t \langle \frac{du(\tau)}{dt},v(\tau)  \rangle_V + \frac{dv(\tau)}{dt},u(\tau)  \rangle_V \, dt.
\end{split}\end{align}
In the Sections \ref{seceasy}, \ref{secapprox} and \ref{secmain} we will 
renounce to mark the integration variables to ensure a better readability.

\subsection{The stress tensor, the couple stress tensor and the electric field}\label{subsectensors}
We denote the symmetric and the skew-symmetric part, resp., of a
tensor $\mathbf{A}$ by $\mathbf{A}^{\sym}:=\frac{1}{2}(\mathbf{A}+\mathbf{A}^\top)$
and $\mathbf{A}^{\anti}:=\frac{1}{2}(\mathbf{A}
-\mathbf{A}^\top)$, respectively. Moreover, we set $\R^{3\times
  3}_{\sym}:=\{\mathbf{A}\in \R^{3\times 3} \fdg
\mathbf{A}=\mathbf{A}^\sym\}$ and $\R^{3\times
  3}_{\anti}:=\{\mathbf{A}\in \R^{3\times 3} \fdg
\mathbf{A}=\mathbf{A}^\anti\}$.  Motivated by the typical example for
the extra stress tensor $\mathbf{S}$ and for the couple stress tensor
$\mathbf{N}$ in (\ref{eq:SN-ex}) and the residual entropy inequality (cf.~\cite[(2.30)]{win-r})
\begin{align}
  \label{eq:cdi}
  \bS:\bD +\bS:\bR + \bN :\nabla \om \ge 0
\end{align}
we make the following assumptions:
\begin{Vss}\label{VssS}
  The stress tensor
  $\mathbf{S}=\mathbf{S}(\mathbf{D},\mathbf{R},\mathbf{E})$ belongs to
  the space \linebreak 
  $C^0(\R_{\sym}^{3\times 3},\R_{\anti}^{3\times 3},\R^3;\setR^{3
    \times 3})$ and fulfills the following assumptions:
  \begin{enumerate}
  \item coercivity:  for all $\bD \in \R_{\sym}^{3\times 3}$,
    $\bR \in \R_{\anti}^{3\times 3}$ and $\bE \in \setR^3$ we have 
   \begin{align}
     \begin{aligned}\label{Skoerziv}
       \mathbf{S}(\mathbf{D},\mathbf{R},\mathbf{E}):\mathbf{D}
       &\geq  c\,\big (1+\vert\E\vert^2\big )\, \big( \vert\mathbf{D}\vert^p-c\big)\,,
       \\
       \mathbf{S}(\mathbf{D},\mathbf{R},\mathbf{E}):\mathbf{R}&\geq
       c\,\vert \mathbf{E}\vert^2 \big( \vert\mathbf{R}\vert^p-c\big)\,,
     \end{aligned}
   \end{align}
 \item boundedness: for all $\bD \in \R_{\sym}^{3\times 3}$,
    $\bR \in \R_{\anti}^{3\times 3}$ and $\bE \in \setR^3$ we have  
   \begin{align}
    \begin{aligned}\label{Sbeschr}
      \vert\mathbf{S}^{\sym}(\mathbf{D},\mathbf{R},\mathbf{E})\vert&\leq
       c\,\big (1+\vert\E\vert^2 \big) \big (1+\vert\mathbf{D}\vert^{p-1}\big)\,,
       \\ 
       \vert\mathbf{S}^{\anti}(\mathbf{D},\mathbf{R},\mathbf{E})\vert&\leq
       c\,\vert \mathbf{E}\vert^2 \big (1+\vert\mathbf{R}\vert^{p-1}\big)\,,
     \end{aligned}
   \end{align}
 \item strict monotonicity: for all $\bD_1, \bD_2 \in
   \R_{\sym}^{3\times 3}$, $\bR_1, \bR_2 \in \R_{\anti}^{3\times3}$
   and $\bE \in \setR^3$ such that $(\mathbf{D}_1,\vert
   \mathbf{E}\vert \mathbf{R}_1)\neq(\mathbf{D}_2,\vert
   \mathbf{E}\vert \mathbf{R}_2)$ we have
   \begin{align}
     \begin{aligned}\label{Smonoton}
       \big (\mathbf{S}(&\mathbf{D}_1,\mathbf{R}_1,\mathbf{E})-
       \mathbf{S}(\mathbf{D}_2,\mathbf{R}_2,\mathbf{E})\big ):
       \big (\mathbf{D}_1-\mathbf{D}_2+\mathbf{R}_1-\mathbf{R}_2\big )>0\,.
     \end{aligned}
   \end{align}
  
 \end{enumerate}
\end{Vss}

\begin{Vss}\label{VssN}
  The couple stress tensor
  $\mathbf{N}=\mathbf{N}(\mathbf{L},\mathbf{E})$ belongs to the space $
  C^0(\R^{3\times3},\R^3;\R^{3\times 3})$ and fulfills the following
  assumptions:
  \begin{enumerate}
  \item coercivity: for all $\bL \in \R^{3\times 3}$ and $\bE \in
    \setR^3$ we have
   \begin{align}\label{Nkoerziv}
      \mathbf{N}(\mathbf{L},\mathbf{E}):\mathbf{L}\geq
      c\,\big
      (1+\vert\mathbf{E}\vert^2
      \big)
      \big(\vert\mathbf{L}\vert^p -c\big)\,,
    \end{align}
  \item boundedness: for all $\bL \in \R^{3\times 3}$ and $\bE \in
    \setR^3$ we have
    \begin{align}\label{Nbeschr}
      \vert\mathbf{N}(\mathbf{L},\mathbf{E})\vert\leq c\,\big
      (1+\vert\mathbf{E}\vert^2
      \big)\big(1+\vert\mathbf{L}\vert^{p-1}\big )\,,
    \end{align}
  \item strict monotonicity: for all $\bL_1,\bL_2 \in \R^{3\times 3}$ and $\bE \in
    \setR^3$ with $\vert \mathbf{E}\vert>0$ and
    $\mathbf{L}_1\neq \mathbf{L}_2$ we have
    \begin{align}\label{Nmonoton}
      (\mathbf{N}(\mathbf{L}_1,\mathbf{E})-\mathbf{N}
      (\mathbf{L}_2,\mathbf{E})):(\mathbf{L}_1-\mathbf{L}_2)>0\,.
    \end{align}
  
  \end{enumerate}
\end{Vss}

\begin{Bem}\label{rem:delta}
{\rm
  We could also have adapted the notion of $(p,\delta)$-structure,
  used e.g.~in \cite{mrs}, \cite{der}, \cite{bdr-7-5} and
  \cite{bdr-phi-stokes}, to the present
  situation. In fact, all results remain valid also under that
  assumption. Moreover, all estimates would depend on $\delta \in
  [\delta_0,\delta_1]$ only through $\delta_0>0$ and $\delta_1$.}
\end{Bem}
In the steady case the quasi-static Maxwell-equations possess very regular solutions. 
Following \cite{rubo}, \cite{frank-phd}, \cite{erw} and the references therein, we know that the 
electric field is a real-analytic function and that and the set $\abs{\gre}^{-1}(0)$ is a finite union of
lower-dimensional $C^1$-manifolds. Especially $\abs{\gre}^{-1}(0)$ is a set of measure zero. If we 
consider the time-dependent case and assume the data to be regular, the solution also possesses 
good regularity properties (cf.~\cite{rubo} for more details). By using Fubini's Theorem we conclude
\begin{align*}
\int_{\Omega_T} \chi_{\gre=\bfzero}= \int_0^T \abs{\gre(t)^{-1}}\, dt =0.
\end{align*}
Therefore we make the following assumption.
\begin{Vss}\label{VssE}
  The electric field $\gre$ belongs to the space
  $L^\infty(0,T;L^\infty(\Omega))$ and a.e.~in $\Omega_T$ there holds
  $\abs{\gre}>0$.
\end{Vss}

Throughout the paper we assume that there exists $p \in (1,\infty)$
and such that $\bS$ satisfies Assumption \ref{VssS} and
$\bN$ satisfies Assumption \ref{VssN}. Moreover, the electric field 
satisfies Assumption \ref{VssE}.

\subsection{Auxiliary results}\label{sec:aux}
In this section we want to present two Lipschitz truncation methods for 
unsteady problems as well as an existence result for parabolic PDEs which 
will be used to solve the easy case $p\geq\frac{11}{5}$ and the approximation 
of our system. The first result is a solenoidal Lipschitz truncation 
which was established in \cite{bds}.

\begin{Sa}[Solenoidal Lipschitz truncation]\label{solpara}
  Let $Q_0=I_0\times B_0$ be a space-time cylinder with a finite
  time-interval $I_0\subset \R$ and a ball $B_0 \subset \R^3$. Let
  $\frac{6}{5}<p<\infty$ and $1 < \sigma <\min\lbrace
  p,p'\rbrace$. Let $(\mathbf{u}_m)$ and $(\mathbf{G}_m)$ satisfy
\begin{align*}
-\int_{Q_0} \mathbf{u}_m\cdot \partial_t \boldsymbol{\xi} \,dx \, dt =
\int_{Q_0} \mathbf{G}_m: \nabla \boldsymbol{\xi}\, dx \, dt \quad \text{for all} 
\ \boldsymbol{\xi} \in C_{0,\Div}^{\infty}(Q_0).
\end{align*}
Assume $(\mathbf{u}_m)$ is a weak null sequence in
$L^p(I_0 ; W^{1,p}(B_0))$, a strong null sequence in $L^\sigma(Q_0)$
and bounded in $L^{\infty}(I_0; L^\sigma(B_0))$.  Further assume that
$\mathbf{G}_m = \mathbf{G}_{1,m}+\mathbf{G}_{2,m}$ is such that
$(\mathbf{G}_{1,m})$ is a weak null sequence in $L^{p'}(Q_0)$ and
$(\mathbf{G}_{2,m})$ converges strongly to zero in $L^\sigma(Q_0)$. Then
there exist double sequences $(\lmk)\subset \R^+$, $(\omk) \subset \R
\times \R^3$, $(\mathbf{u}_{m,k})\subset L^1(Q_0)$ and $k_0\in \Nat$
such that for all $k\geq k_0$: 
\begin{enumerate}
\item[\textnormal{(a)}] $2^{2^k} \leq  \lmk \leq 2^{2^{k+1}}$.
\item[\textnormal{(b)}]
  $(\mathbf{u}_{m,k}) \subset L^s(\frac{1}{4}I_0 ;
  W_{0,\Div}^{1,s}(\frac{1}{6}B_0 ))$
  for all $s<\infty$ and
  $\operatorname{supp}(\mathbf{u}_{m,k}) \subset
  \frac{1}{6}Q_0$.\footnote{Let
    $I_0=(t_0-\rho,t_0+\rho)$ and $B_0=B_r(x_0)$. Then we define the
    scaled space-time cylinder $\frac{1}{6}Q_0$ by
    $(t_0-\frac{1}{6}\rho,t_0+\frac{1}{6}\rho)\times
    B_{\frac{1}{6}r}(x_0)$. }
\item[\textnormal{(c)}] $\mathbf{u}_{m,k} = \mathbf{u}_{m}$ a.e.~on $\frac{1}{8}Q_0 \setminus \omk$.
\item[\textnormal{(d)}] $\norm{\nabla \mathbf{u}_{m,k}}_{L^\infty(\frac{1}{4}Q_0)} \leq c\, \lmk$.
\item[\textnormal{(e)}] $ \mathbf{u}_{m,k} \rightarrow \bfzero$ in $ L^\infty(\frac{1}{4}Q_0)$ 
for $m \rightarrow \infty $ and $k$ fixed.
\item[\textnormal{(f)}] $\nabla \mathbf{u}_{m,k} \rightharpoonup \bfzero$ in $ L^s(\frac{1}{4}Q_0)$ 
for $m \rightarrow \infty $ and $k$ fixed.
\item[\textnormal{(g)}] $\limsup_{m\rightarrow \infty} \lmk^p \abs{\omk} \leq c\,2^{-k}$.
\item[\textnormal{(h)}]
  $\limsup_{m\rightarrow \infty} \abs{ \int \mathbf{G}_m : \nabla
    \mathbf{u}_{m,k}\, dx \, dt } \leq c \,\lmk^p \abs{\omk}$.
\end{enumerate}
\end{Sa}
\begin{Bew}
This is Theorem 2.2 in \cite{bds}. Another proof can be found in \cite{erik-dipl} 
where the result of \cite[Lemma 2.6]{bds} is proofed differently.
\end{Bew}
We also cite a useful corollary of this theorem, cf.~\cite[Corollary 2.4]{bds}.
\begin{Kor}\label{corsolpara}
  Let all the assumptions of Theorem \ref{solpara} be satisfied and
  let $\zeta \in C_0^{\infty}(\frac{1}{6}Q_0)$ satisfy
  $\chi_{\frac{1}{8}Q_0}\leq \zeta \leq \chi_{\frac{1}{6}Q_0} $. Then
  for every $\vek{K} \in L^{p'}(\frac{1}{6}Q_0)$
\begin{align*}
\limsup\limits_{m\rightarrow \infty} \left\vert\int ((\mathbf{G}_{1,m}+\vek{K}):\nabla 
\mathbf{u}_m)\,\zeta\, \chi_{\omk^\complement} \, dx \, dt \right\vert \leq c\, 2^{-\frac{k}{p}}.
\end{align*}
\end{Kor}
In \cite{die-ru-wolf} a Lipschitz truncation method for non-solenoidal
functions was developed.  Since we need a small generalization of
\cite[Theorem 3.9]{die-ru-wolf} we sketch the proof.  Let us start
with some notation.  For $\alpha >0$ we define the anisotropic
parabolic metric $d_{\alpha}$ in $\R \times \R^3$ by
$d_{\alpha}((t,x),(s,y)):= \max \lbrace \abs{x-y}, \abs{\alpha^{-1}
  (t-s)}^{\frac{1}{2}} \rbrace$.
For $(t,x) \in \R \times \R^3$ and $r>0$ we define $\alpha$-parabolic
cylinders
$Q_r^\alpha((t,x)):= \lbrace (s,y) \in \R \times \R^3 \fdg
d_{\alpha}((t,x),(s,y))<r\rbrace$.
Note that
$ Q_r^\alpha((t,x))= (t-\alpha r^2,\linebreak t+ \alpha r^2) \times B_r(x)$.  Let
$E \subset \R \times \R^3$ be open and bounded, where $\R \times \R^3$
is equipped with the anisotropic metric $d_{\alpha}$ for some
$\alpha>0$.  According to \mbox{\cite[Lemma 3.1]{die-ru-wolf}} there exists
a Whitney type covering $(Q_i^\alpha)_{i \in \Nat}$ of $E$ with
$\alpha$-parabolic cylinders $Q_i^\alpha:= Q_{r_i}^\alpha((t_i,x_i))$.
There also exists a subordinate partition of unity
$(\psi_i)_{i \in \Nat}$ to this Whitney type covering
$(Q_i^\alpha)_{i \in \Nat}$ (cf.~\cite[(3.2)]{die-ru-wolf}).  So we
are able to define the truncation operator.
\begin{Def}\label{truncdef}
Let $E \subset  \Omega_T$ be open. For 
$\vek{u} \in L^1_{\operatorname{loc}} (\Omega)T)$ we define
\begin{align*}
(\trunc\vek{u}) ((t,x)) := \begin{cases} 
u(t,x) \qquad \qquad &\textnormal{if } (t,x) \in \Omega_T \setminus E   
\\ \sum\limits_{i \in \Nat} \psi_i(t,x) u_{Q_i^{\alpha}} &\textnormal{if } (t,x) \in E,
\end{cases}
\end{align*}
where $u_{Q_i^{\alpha}}:= \dint_{Q_i^{\alpha}} u \, dx \, dt $ is the mean value over $Q_i^{\alpha}$.
\end{Def}
With this definition we get the first continuity result.
\begin{Lem}\label{truncstetig}
There exists a constant $c$ such that for every $1\leq p \leq  \infty$
\begin{align*}
\norm{\trunc u}_{L^p( \Omega_T)} \leq c\, \norm{u}_{L^p( \Omega_T)}.
\end{align*}
\end{Lem}
\begin{Bew}
cf.~\cite[Lemma 3.5]{die-ru-wolf}.
\end{Bew}
Before we state the main theorem of chapter 3 in \cite{die-ru-wolf} in a generalized version, 
we need some results for the maximal operator. For $g \in L^p(\R^{1+3}) \ (p>1)$  we denote by
\begin{align*}
\mathcal{M}_x(f)(t,x)&:= \sup\limits_{0<r<\infty} \dint_{B_r(x)} \abs{f(t,y)}dy, \\ 
\mathcal{M}_t(f)(t,x)&:= \sup\limits_{0<\rho<\infty} \dint_{I_\rho(t)} \abs{f(s,x)}ds,
\end{align*}
the usual maximal operators in the space and time variables. Here 
$I_\rho(t):= (t-\rho,t+\rho)$. More details concerning maximal operators 
can be found in \cite{stein}. We want to use the composition of
these two maximal operators
\begin{align*}
\mathcal{M}^\ast(g):= \mathcal{M}_t(\mathcal{M}_x(f)).
\end{align*}
This maximal operator satisfies strong and weak type estimates and for 
all $(t,x) \in \R^{1,3}$ and $r,\rho>0$ there holds (cf.~\cite[Appendix A]{die-ru-wolf})
\begin{align}\label{g2.3max}
\dint_{I_\rho(t)} \dint_{B_r(x)} \abs{g(s,y)}\,dy\,ds \leq \mathcal{M}^\ast(g)(t,x).
\end{align}
\begin{Sa}\label{paralip}
For $1<p,\sigma<\infty$ let $\mathbf{u} \in L^\infty(I;L^2(\Omega))\cap L^p(I;W^{1,p}(\Omega))$, 
$\mathbf{H} \in L^\sigma(\Omega_T)$, $\mathbf{k} \in L^{p'}(\Omega_T)$ satisfy 
\begin{align}\label{g2.3.3}
-\int_{\Omega_T} \mathbf{u}\cdot \partial_t \tphi \, dx \,dt=\int_{\Omega_T} \mathbf{H}: 
\nabla \tphi \, dx \,dt + \int_{ \Omega_T} \mathbf{k}\cdot\tphi \, dx \,dt
\end{align}
for all $\tphi \in C_0^\infty(\Omega_T)$.
We define $(\Lambda>0)$
\begin{align*}
\mathcal{O}_\Lambda&:= \lbrace (t,x) \in \R^{1+3} \fdg \mathcal{M}^{\ast}(\abs{\nabla \mathbf{u}})(t,x)+
\alpha\mathcal{M}^{\ast}(\abs{\mathbf{H}})(t,x)+
\alpha\mathcal{M}^{\ast}(\abs{\mathbf{k}})(t,x)>\Lambda\rbrace,\\
\mathcal{U}_1&:=\lbrace (t,x) \in \R^{1+3} \fdg  \mathcal{M}^{\ast}(\abs{\mathbf{u}})(t,x)>1\rbrace.
\end{align*}
Let $E$ be an open, bounded set such that $\Omega_T\cap(\mathcal{O}_\Lambda\cup \mathcal{U}_1)\subset E\subset \Omega_T$.
Let $K \subset \Omega_T$ be a compact set, then there holds:
\begin{compactenum}
\item[\textnormal{(i)}]\label{paralip1} The Lipschitz truncation $\trunc\mathbf{u}$ belongs to $C_{d_{\alpha}}^{0,1}(K)$, 
the space of Lipschitz continous functions with respect to $d_\alpha$ on $K$, 
where the norm depends on $K$, $\Lambda$, $\alpha$, $\norm{\mathbf{u}}_{L^1(E)}$ 
and $\norm{\mathbf{u}}_{L^1(\Omega_T)}$. 
In particular $\trunc\mathbf{u}, \nabla\trunc\mathbf{u} \in L^\infty(K)$.
\item[\textnormal{(ii)}]\label{paralip2}  The Lipschitz truncation $\trunc \mathbf{u}$ satisfies the estimates
\begin{align*}\begin{split}
\norm{\nabla \trunc \mathbf{u}}_{L^\infty(K)}&\leq
 c\,(\Lambda+\alpha^{-1}\delta_{\alpha,K}^{-3-3}\norm{\mathbf{u}}_{L^{1}(E)}),\\
\norm{\trunc \mathbf{u}}_{L^\infty(K)}&\leq
 c\,(1+\alpha^{-1}\delta_{\alpha,K}^{-3-2}\norm{\mathbf{u}}_{L^{1}(E)}),
\end{split}\end{align*}
where $\delta_{\alpha,K}:=d_{\alpha}(K,\partial  \Omega_T)$.
\item[\textnormal{(iii)}]\label{paralip3}  The function $(\partial_t\trunc\mathbf{u})\cdot( \trunc\mathbf{u} -\mathbf{u} )$
belongs to $L^1(K\cap E)$ and we have
\begin{align*}\begin{split}
\norm{(\partial_t\trunc\mathbf{u})\cdot( \trunc\mathbf{u} -\mathbf{u} )}_{L^1(K\cap E)} \leq
 c\,\alpha^{-1}\abs{E}(\Lambda+\alpha^{-1}  \delta_{\alpha,K}^{-3-3} \norm{\mathbf{u}}_{L^{1}(E)} )^2\,.
\end{split}\end{align*}
\item[\textnormal{(iv)}]\label{paralip4}  For all $\zeta \in C_0^\infty(\Omega_T)$ there holds the identity 
\begin{align*}\begin{split}
&\int_I \left\langle \frac{d\mathbf{u}}{dt}(t), ( \trunc\mathbf{u}(t) )\zeta(t)  \right\rangle dt\\
&=\frac{1}{2} \int_{ \Omega_T} ( \abs{\trunc\mathbf{u}}^2 -2 \mathbf{u}\cdot\trunc\mathbf{u}  ) \partial_t \zeta \, dx\, dt + \int_E (\partial_t\trunc\mathbf{u}    )\cdot (\trunc\mathbf{u} -\mathbf{u}   )\zeta \,dx\, dt\,,
\end{split}\end{align*}
where $\langle \cdot,\cdot \rangle$ denotes the usual duality pairing with respect to $\Omega$.
\end{compactenum}
\end{Sa}
In \cite{die-ru-wolf} this Theorem is proved only with $\mathbf{k}\equiv \bfzero$. 
In order to deal with $\mathbf{k} \neq \bfzero$ we changed the definition of $\mathcal{O}_{\Lambda}$, which  
coincides with the definition in \cite{die-ru-wolf} for $\mathbf{k}\equiv \bfzero$.
The proof of \cite[Theorem 3.9]{die-ru-wolf} uses a Poincar\'{e}-type
inequality (cf.~\cite[Appendix B]{die-ru-wolf}, \cite[Lemma B.3]{wolf-diss}) 
which allows to control mean values of the form\footnote{Here for $Q_r^\alpha((t,x))=\lbrace (s,y) \fdg d_\alpha( (t,x), (s,y)  ) <r \rbrace$ we define the scaled $\alpha$-parabolic cylinder by  $4Q_r^\alpha((t,x))=\lbrace (s,y) \fdg d_\alpha( (t,x), (s,y)  ) <4r \rbrace$.} 
$\dint_{Q_i} \abs{\mathbf{u} - \mathbf{u}_{Q_i}}\, dx\, dt$ by norms of
$\nabla \mathbf{u}$, $\mathbf{H}$ and $\mathbf{k}$.
Thus, if we generalize the Poincar\'{e}-type inequality for functions 
with a distributional time derivative of the form of \eqref{g2.3.3}, 
the proof of \cite[Theorem 3.9]{die-ru-wolf} can be applied to prove Theorem \ref{paralip}.
\begin{Lem}[Poincar\'{e}-type inequality]\label{pointype} For $\alpha,r>0$ let
 $Q_r^\alpha:=Q_r^\alpha((t,x))$ be a $\alpha$-parabolic cylinder.
Assume $\bfu \in L^1(Q_r^\alpha)$ with $\nabla \bfu \in L^1(Q_r^\alpha)$, $\bfH \in L^1(Q_r^\alpha)$ 
and $\bfk \in  L^1(Q_r^\alpha)$ satisfy
\begin{align}\label{g2.3.4}
-\int_{Q_r^\alpha} \mathbf{u}\cdot \partial_t \tphi \, dx \,dt=\int_{Q_r^\alpha} \mathbf{H}: \nabla \tphi \, dx \,dt + \int_{Q_r^\alpha} \mathbf{k}\cdot\tphi \, dx \,dt
\end{align}
for all $\tphi \in C_0^{\infty}(Q_r^\alpha)$. Then there exists a constant $c$ independent of $\alpha$, $r$ and $(t,x)$ such that
\begin{align*}
\int_{Q_r^\alpha} \abs{\mathbf{u} - \mathbf{u}_{Q_r^\alpha}}\, dx \,dt \leq 
c\,r\left( \norm{\nabla\bfu}_{L^1(Q_r^\alpha)}+ \alpha \,\norm{\bfH}_{L^1(Q_r^\alpha)} 
 + \alpha \,r \,\norm{\bfk}_{L^1(Q_r^\alpha)} \right).
\end{align*}   
\end{Lem}
\begin{Bew}
We only proof the special case $(t,x)=(0,0)$, $1=\alpha=r$ since the general result of 
this lemma follows from a transformation of coordinates. We abbreviate 
$Q_1:=Q_1^1((0,0))$ and $B_1:= B_1(0)$. With the same arguments as in 
\cite{die-ru-wolf} we get 
\begin{align}\label{g2.3.5}
\int_{Q_1}\! \abs{\bfu(t,x)- \bfu_{Q_1}} \, dx \,dt \leq c \int_{Q_1}\! \abs{\nabla \bfu}  \, dx \,dt +
\frac{\abs{B_1}}{2} \int_{-1}^1\int_{-1}^1\! \big\vert  \bfg \big( \bfu(t) -\bfu(s)\big) \big\vert \, ds \,dt  
\end{align}
where $\bfg:L^1(B_1) \rightarrow \R^3$ is defined for $\zeta \in C_0^\infty(B_1)$ satisfying $\chi_{\frac{1}{2}B_1} \leq \zeta\leq  \chi_{B_1}$, by 
\begin{align*}
\bfg(\bfv) =\frac{1}{\int_{B_1} \zeta \, dx}\int_{B_1} \zeta \bfv \, dx.
\end{align*}
For arbitrary $\x \in C_0^\infty(B_1)$ and $\gamma \in  C_0^\infty(-1,1)$ the function 
$\tphi(t,x):= \x(x)\gamma_{-h}(t)$, where $\gamma_{-h}(t):= \frac{1}{h}\int_t^{t-h}\gamma(s) \, ds$ 
is the Steklov average, is a admissible testfunction for equation \eqref{g2.3.4}. 
By standard arguments concerning the Steklov average we conclude
\begin{align*}
\int_{-1}^1 \int_{B_1} \partial_t \bfu_h\cdot \x \, dx \,\gamma\, dt= \int_{-1}^1 \int_{B_1} \bfH_h:\nabla\x\, dx \,\gamma\, dt +\int_{-1}^1 \int_{B_1} \bfk_h\cdot \x\, dx\, \gamma \, dt
\end{align*}
and with the fundamental lemma of the calculus of variations we get
\begin{align*}
\int_{B_1} \partial_t \bfu_h(\tau)\cdot \x \, dx=  \int_{B_1} \bfH_h(\tau):\nabla\x\, dx +\int_{B_1} \bfk_h(\tau)\cdot \x\, dx
\end{align*}
for almost all $\tau \in (-1,1)$.
Splitting this equation into the components $\bfu_h^i$, $\bfk_h^i$ and $\bfH_h^i$, 
 where $\bfH_h^i$ is the $i$-th line of the tensor $\bfH_h$, 
 by using a testfunction which is 0 in all components except the
 $i$-th component, we get
\begin{align*}
\int_{B_1} \partial_t \bfu_h^i(\tau)\, \xi \, dx=  \int_{B_1} \bfH_h^i(\tau)\cdot\nabla\xi\, dx +\int_{B_1} \bfk_h^i(\tau)\, \xi\, dx \quad \text{for all } \xi \in C_0^\infty(B_1).
\end{align*}
This equation, the properties of the Steklov average, the definition of $\bfg$ and 
the properties of $\zeta$ imply
\begin{align*}\begin{split}
\big\vert \bfg\big( \bfu(t) -\bfu(s)\big)\big\vert&=\lim\limits_{h\rightarrow 0} \big\vert \bfg\big( \bfu_h(t) -\bfu_h(s)\big)\big\vert\\
&= \lim\limits_{h\rightarrow 0} \Big\vert \frac{1}{\int_{B_1} \zeta \, dx}\int_s^t\int_{B_1} \partial_t \bfu_h \zeta  \, dx \, d\tau \Big\vert\\
&\leq \lim\limits_{h\rightarrow 0} \frac{c}{\int_{B_1} \zeta \, dx} \sum_{i=1}^3 \Big\vert\int_s^t\int_{B_1} \partial_t \bfu_h^i \zeta  \, dx \, d\tau \Big\vert\\
&\leq \lim\limits_{h\rightarrow 0} \frac{c}{\abs{\frac{1}{2}B_1}} \sum_{i=1}^3 \Big\vert\int_s^t\int_{B_1}  \bfH_h^i:\nabla \zeta +  \bfk_h^i \zeta \, dx \, d\tau \Big\vert\\
&\leq \lim\limits_{h\rightarrow 0}\norm{\bfH_h}_{L^1(Q_1)}+ \norm{\bfk_h}_{L^1(Q_1)} \\
&= \norm{\bfH}_{L^1(Q_1)}+ \norm{\bfk}_{L^1(Q_1)}
\end{split}\end{align*}
This together with \eqref{g2.3.5} proofs the special case of this lemma.
\end{Bew}
With this Poincar\'{e}-type inequality we can prove \cite[Lemma 3.11]{die-ru-wolf} which is 
important for the proof of \cite[Theorem 3.9]{die-ru-wolf}
\begin{Lem}\label{DRWLem3.11.}
Under the assumptions of Theorem \ref{paralip} we have for all $Q_i^\alpha$ belonging to the
 Whitney covering of $E$ such that $Q_i^\alpha \cap K \neq \emptyset $
\begin{align}
  \dint_{4Q_i^\alpha}\abs{\bfu - \bfu_{4Q_i^\alpha}}\,dx\,dt \leq c\,
  r_i\big (\Lambda +\alpha^{-1} \delta_{\alpha,K}^{-3-3}\norm{\bfu}_{L^1(E)}\big),
\end{align}
where the constant depends on the diameter of $\Omega$.
\end{Lem}
\begin{Bew}
The properties of the Whitney covering imply (i) $16Q_i^\alpha\cap \big( (\mathcal{O}_{\Lambda})^\mathsf{c} 
\cap (\mathcal{U}_1)^\mathsf{c} \big)$ or (ii) $16Q_i^\alpha\cap Q\neq \emptyset$. In case (i) we use the new 
Poincar\'{e}-type inequality in Lemma \ref{pointype} to estimate 
\begin{align*}
\dint_{4Q_i^\alpha}\abs{\bfu - \bfu_{4Q_i^\alpha}}\,dx\,dt \leq c
  \,r_i \dint_{4Q_i^\alpha} \abs{\nabla \bfu} +\alpha\,
  \abs{\bfH}+\alpha \,r_i\abs{\bfk}\,dx\,dt, 
\end{align*}
where we used $4Q_i^\alpha \subset E \subset Q$, which also implies 
$0\leq r_i \leq \operatorname{diam}\Omega$. Since 
$16Q_i^\alpha\cap (\mathcal{O}_{\Lambda})^\mathsf{c} \neq \emptyset$ we 
find a $(\widetilde{t},\widetilde{x})\in 16Q_i^\alpha\cap (\mathcal{O}_{\Lambda})^\mathsf{c}$. 
The definition of the $\alpha$-parabolic cylinders imply 
$4Q_i^\alpha\subset Q_{20r_i}^{\alpha}(t_0,x_0) $. Now from the above inequality, 
the definition of the maximal operator $\mathcal{M}^\ast$, \eqref{g2.3max} and 
the new definition of $\mathcal{O}_\Lambda$ in Theorem \ref{paralip} we get
\begin{align*}
  \dint_{4Q_i^\alpha}\abs{\bfu - \bfu_{4Q_i^\alpha}}\,dx
  &\,dt \leq c\,r_i \,\dint_{4Q_i^\alpha} \abs{\nabla \bfu}
    +\alpha\, \abs{\bfH}+\alpha \, r_i\abs{\bfk}\,dx\,dt
  \\
  &\leq c(\Omega)\, r_i
    \,\dint_{Q_{20r_i}^{\alpha}(\widetilde{t},\widetilde{x})}
    \abs{\nabla \bfu} +\alpha \,\abs{\bfH}+\alpha\, \abs{\bfk}\,dx\,dt
 \\
  &\leq c(\Omega)\, r_i \big(\mathcal{M}^{\ast}(\nabla\bfu)(\widetilde{t},\widetilde{x})
    +\alpha\mathcal{M}^{\ast}(\bfH)(\widetilde{t},\widetilde{x})+
    \alpha\mathcal{M}^{\ast}(\bfk)(\widetilde{t},\widetilde{x})\big)
  \\ 
  &\leq c(\Omega)\, r_i\, \Lambda.
\end{align*}
Case (ii) can be treated exactly as in \cite{die-ru-wolf}.
\end{Bew}
Now the proof of Theorem \ref{paralip} is exactly the same as the proof of \cite[Theorem 3.9]{die-ru-wolf} 
we just have to use Lemma \ref{DRWLem3.11.} whenever \cite[Lemma 3.11]{die-ru-wolf} is used.

\smallskip
Finally we quote an existence result for parabolic PDEs (cf.~\cite{br-hirano}).
\begin{Sa}\label{ResultatHirano}
Let $(V,H,V^\ast)$ be a Gelfand-Triple. Assume $Z$ is another reflexive, separable Banach space such that
$Z\hookrightarrow V$ with a continuous and dense embedding. Moreover, assume that there exists an increasing
sequence of finite dimensional subspaces $V_n \subseteq Z$, such that 
$\cup _{n \in \mathbb{N}} V_n$ is dense in $V$. 
Additionally, there exists self-adjoint projections $P_n:H \rightarrow H$, such that $P_n(V)=V_n$ and
$\norm{P_{n\ |Z}}_{\mathcal{L}(Z,Z)}\leq c$ with a constant $c$ independent of $n\in\mathbb{N}$.
Finally, let $\lbrace A(t) \fdg 0\leq t\leq T\rbrace$ be a family of operators 
from $V$ to $V^\ast$ with the following properties:
\begin{itemize}
\item[\textnormal{(A1)}] $A(t):V \rightarrow V^\ast$ is pseudomonotone for almost every $t \in [0,T]$.
\item[\textnormal{(A2)}] For every $u \in L^p(0,T;V)\cap L^\infty(0,T;H)$ the mapping $t \mapsto A(t)u(t)$ 
from $[0,T]$ to $V^\ast$ is Bochner-measurable.
\item[\textnormal{(A3)}] There exists a positive constant $c_1$ and a nonnegative 
function $C_2 \in L^1(0,T)$,
such that
\begin{align*}
\langle A(t)x,x  \rangle_V \geq c_1 \norm{x}_V^{p} -C_2(t)
\end{align*}
for almost every $t \in [0,T]$ and all $x\in V$.
\item[\textnormal{(A4)}] There exists $0< q<\infty$, as well as constants $c_3>0$, $c_4\geq 0$ 
and a nonnegative function $C_5\in L^{p'}(0,T)$,
such that
\begin{align*}
\norm{A(t)x}_{V^\ast}\leq c_3\norm{x}_V^{p-1}+ c_4 \norm{x}_H^q \norm{x}_V^{p-1} +  C_5(t)
\end{align*}
for almost every $t \in [0,T]$ and all $x\in V$.
\end{itemize}
Then for every $u_0 \in H$, $f \in L^{p'}(0,T;V^\ast)$ there exists a function 
$u \in W_p^1(0,T;V,H)$ such that
\begin{align*}
u'(t) +A(t)u(t)&=f(t) &&\hspace*{-15mm}\text{in }V^\ast \text{ for a.e.} \ t \in [0,T],\\
u(0)&=u_0 &&\hspace*{-15mm}\text{in}\ H.
\end{align*}
\end{Sa}
\section{Easy case $p\geq \frac{11}{5}$} \label{seceasy}

\begin{Sa}\label{easycase}
  Let $p\in [\frac{11}{5},\infty)$, $T \in (0,\infty)$ and
  $\Omega\subset\R^3$ be a bounded domain with
  Lipschitz-boundary. Assume that $\grs$ satisfies Assumption
  \ref{VssS}, that $\grn $ satisfies Assumption \ref{VssN} and that $ \gre$
  satisfying Assumption \ref{VssE} is given.  Then there exists for all
  $\vvv_0\in H(\Omega)$, $\om_0 \in L^2(\Omega)$ and
  $\f,\l \in L^{p'}(\Omega_T)$ a weak solution
  $(\vvv, \om) \in W_p^1(0,T;V_p(\Omega),H(\Omega)) \times
  W_p^1(0,T;W_0^{1,p}(\Omega),L^2(\Omega))$
  of the problem \eqref{system}--\eqref{ic} satisfying for all
  $(\tphi, \tpsi) \in L^p(0,T;V_p(\Omega))\times
  L^p(0,T;W_0^{1,p}(\Omega))$
  \begin{align}
    \begin{split}\label{easycaseeq}
      &\int_0^T \langle \frac{d\vvv}{dt}(t),\tphi(t) \rangle_{V_p}+\int_{\omt}\sdre:\nabla\tphi-\int_{\omt}\vvv\otimes\vvv:\nabla\tphi\\
      &\quad+\int_0^T \langle \frac{d\om}{dt}(t),\tpsi(t) \rangle_{W_0^{1,p}} +\int_{\omt} \grn(\nabla\om,\gre):\nabla\tpsi- \int_{\omt}\om\otimes\vvv:\nabla\tpsi\\
      &=\int_{\omt} \f\cdot\tphi +\int_{\omt} \l \cdot
      \tpsi-\int_{\omt} \big (\levy:\sdre \big)\cdot \tpsi
    \end{split}
  \end{align}
  as well as $\vvv(0)=\vvv_0$ and $\om(0)=\om_0$.
\end{Sa}
\begin{Bew}
We want to use Theorem \ref{ResultatHirano}. In view of
\eqref{eq:cdi}, the identity \eqref{rulelevy}
and the assumptions on the stress tensors it is natural to view the
system \eqref{system} as a unit. Thus we are searching two unknown 
functions $\vvv$ and $\om$ as elements $(\bv , \om)$ of the product
space $V_p(\Omega) \times W_0^{1,p}(\Omega)$. To simplify the notation
we set  
\begin{alignat*}{2}
  \begin{aligned}
    \mathscr{V}_p&:=V_p(\Omega) \times W_0^{1,p}(\Omega), \quad
    &&\norm{(\vek{u},\vek{w})}_{\mathscr{V}_p} := \big(
    \norm{\vek{u}}_{V_p}^2 + \norm{\vek{w}}_{W_0^{1,p}}^2
    \big)^\frac{1}{2},
    \\ 
    \mathscr{H}&:= H(\Omega) \times L^2(\Omega),
    &&\norm{(\vek{u},\vek{w})}_{\mathscr{H}} := \big(
    \norm{\vek{u}}_{H}^2 + \norm{\vek{w}}_{L^2}^2 \big)^\frac{1}{2}.
  \end{aligned}
\end{alignat*}
Since $(V_p(\Omega), H(\Omega),V_p(\Omega)^\ast)$ and
$(W_0^{1,p}(\Omega), L^2(\Omega),W_0^{1,p}(\Omega)^\ast)$ form
Gelfand-trip\-les, it is obvious that
$(\mathscr{V}_p,\mathscr{H},(\mathscr{V}_p)^\ast)$ forms a
Gelfand-triple as well. We set
$\mathscr{Z}:= V^3(\Omega)\times W_0^{3,2}(\Omega)$. Then, according
to \cite[Appendix 4.11 and 4.14]{MNRR}, we know that
$\mathscr{Z},\mathscr{V}_p, \mathscr{H}$
satisfy all assumptions in Theorem \ref{ResultatHirano}. \\
Next we define operators
$A(t),A_i(t):\mathscr{V}_p \rightarrow(\mathscr{V}_p)^\ast$,
$i=1,\ldots,4$, by
\begin{align}\begin{split}\label{defOperatoren}
\left\langle A_1(t)(\vek{u},\vek{w}),(\tphi,\tpsi) \right\rangle&:= 
\int_\Omega \grs\big(\vek{D}\vek{u},\grr(\vek{u},\vek{w}),\gre(t)\big) :(\vek{D}\tphi + \grr(\tphi,\tpsi)), \\
\left\langle A_2(t)(\vek{u},\vek{w}),(\tphi,\tpsi) \right\rangle&:= - \int_\Omega \vek{u}\otimes\vek{u}:\nabla\tphi, \\
\left\langle A_3(t)(\vek{u},\vek{w}),(\tphi,\tpsi) \right\rangle&:=\int_\Omega \grn(\nabla\vek{w},\gre(t)):\nabla\tpsi,   \\
\left\langle A_4(t)(\vek{u},\vek{w}),(\tphi,\tpsi) \right\rangle&:=-\int_\Omega \vek{w}\otimes\vek{u}:\nabla\tpsi,  \\
A(t)&:= A_1(t)+A_2(t)+A_3(t)+A_4(t),
\end{split}\end{align}
where $\langle \cdot ,\cdot \rangle$ denotes the duality pairing between $\mathscr{V}_p$ and $(\mathscr{V}_p)^\ast$.
In order to apply Theorem \ref{ResultatHirano} we have to check that the family $A(t)$ 
satisfies  (A1)--(A4).

(A1): Using \eqref{Sbeschr}, \eqref{Nbeschr} as well as
$E(t)\in L^\infty(\Omega)$ for almost every $t \in [0,T]$ we conclude
with the theory of Nemyckii operators that $A_1(t)$ and $A_3(t)$ are
continuous operators. Moreover, from \eqref{Smonoton} and
\eqref{Nmonoton} we get the monotonicity of $A_1(t)$ and
$A_3(t)$. Using the compact embedding
$\mathscr{V}_p\hookrightarrow\hookrightarrow L^s(\Omega) \times
L^s(\Omega)$
for any $1\leq s <\frac{3p}{3-p}$ and H{\"o}lder's inequality it is
easy to show that $A_2(t)$ and $A_4(t)$ are strongly continuous
operators for $p >\frac{9}{5}$. Since monotone, continuous operators
and strongly continuous operators are pseudomonotone and since
summation maintains pseudomonotonicity, we conclude that $A(t)$ is
pseudomonotone for almost every $t \in [0,T]$.

(A2): Fix $(\vek{u},\vek{w})\in L^p(0,T;\mathscr{V}_p)\cap L^\infty(0,T;\mathscr{H})$. The 
argumentation will be the same as in \cite{br-hirano}. Since $\mathscr{V}_p$ is separable we are
able to use Pettis' Theorem. Therefore it is enough to prove that 
$t\mapsto \left\langle A(t)(\vek{u}(t),\vek{w}(t)),(\tphi,\tpsi) \right\rangle_{\mathscr{V}_p}$
is Lebesgue measurable for arbitrary $(\tphi,\tpsi)\in \mathscr{V}_p$. Using 
\eqref{Sbeschr}, \eqref{Nbeschr}, \eqref{VssE} it is clear that every function 
appearing in the definitions of $A_i$ is an element of $L^1(\Omega_T)$ so that 
Fubini's theorem yields the assertion.

(A3) For arbitrary $(\vek{u},\vek{w}) \in \mathscr{V}_p$, $p \ge\frac 95$, there holds
\begin{align}\begin{split}\label{g3.0a}
\left\langle A_2(t)(\vek{u},\vek{w}),\vek{u},\vek{w}) \right\rangle &= - \int_\Omega \vek{u}\otimes\vek{u}:\vek{u}=0,\\
\left\langle A_4(t)(\vek{u},\vek{w}),\vek{u},\vek{w}) \right\rangle&=-\int_\Omega \vek{w}\otimes\vek{u}:\vek{w}=0 , 
\end{split}\end{align}
due to $\operatorname{div}\vek{u}=0$ and integration by parts. Using \eqref{Skoerziv} and \eqref{Nkoerziv} 
we immediately estimate
\begin{align}\begin{split}\label{g3.0}
&\left\langle A_1(t)(\vek{u},\vek{w}),(\vek{u},\vek{w}) \right\rangle+\left\langle A_3(t)(\vek{u},\vek{w}),(\vek{u},\vek{w}) \right\rangle\\
&\geq c\int_{\Omega} \abs{\grd\vek{u}}^p  +  c\int_{\Omega} \abs{\nabla\vek{w}}^p - c(\gre,\Omega)\\
&= c \norm{\vek{u}}_{V_p}^p + c \norm{\vek{w}}_{W_0^{1,p}}^p -c(\gre,\Omega)\geq c \norm{(\vek{u},\vek{w})}_{\mathscr{V}_p}^p-c(\gre,\Omega),
\end{split}\end{align}
where the constant $c(\gre,\Omega)$ depends only on
$\norm{\bE}_{L^\infty(\Omega_T)}$ and $\abs{\Omega}$ because of
Assumption \ref{VssE}.

(A4) Let us start with the operator $A_1$. Using \eqref{Sbeschr}, H{\"o}lder's inequality and the 
continuous embedding $L^p\hookrightarrow L^1$ we get
\begin{align}\begin{split}\label{g3.1}
\norm{A_1(t) (\vek{u},\vek{w})}_{(\mathscr{V}_p)^\ast}
&\leq \sup\limits_{\norm{(\tphi,\tpsi)}_{\mathscr{V}_p}\leq1  }
\int_\Omega c\,(1+\abs{\gre(t)}^2)(1+ \abs{\grd\vek{u}}^{p-1}) \abs{\grd\tphi}\\
&\qquad+\sup\limits_{\norm{(\tphi,\tpsi)}_{\mathscr{V}_p}\leq1  }\int_\Omega c\,\abs{\gre(t)}^2(1+ \abs{\grr(\vek{u},\vek{w})}^{p-1}) \abs{\grr(\tphi,\tpsi)}\\
&\leq \sup\limits_{\norm{(\tphi,\tpsi)}_{\mathscr{V}_p}\leq1  } 
c_{\gre,p,\abs{\Omega}}(1+\norm{\grd\vek{u}}_{L^p}^{p-1})\norm{\grd\tphi}_{L^p}\\
&\qquad+\sup\limits_{\norm{(\tphi,\tpsi)}_{\mathscr{V}_p}\leq1  }
c_{\gre,p,\abs{\Omega}}(1+\norm{\grr(\vek{u},\vek{w})}_{L^p}^{p-1})\norm{\grr(\tphi,\tpsi)}_{L^p}.\raisetag{39mm}
\end{split}\end{align}
Here the constant $c_{\gre,p,\abs{\Omega}}$ is again independent of $t\in [0,T]$ because of Assumption \ref{VssE}.
Due to \eqref{levyabsch} we are able to estimate that for any $(\tphi,\tpsi) \in \mathscr{V}_p$ and 
 a.e. $x \in \Omega$
\[ \abs{ \grr(\tphi(x),\tpsi(x))} \leq \abs{\nabla\tphi(x)} +c \abs{\tpsi(x)} .     \]
This, together with Poincare's and Korn's inequality, implies
\begin{align}\begin{split}\label{g3.2}
\norm{\grr(\tphi,\tpsi)}_{L^p} &\leq c\, ( \norm{\grd \tphi}_{L^p}+ \norm{\nabla\tpsi}_{L^p})
= c\, ( \norm{\tphi}_{V_p}+ \norm{\tpsi}_{W_0^{1,p}})\\
&\leq c\,\norm{(\tphi,\tpsi)}_{\mathscr{V}_p}.
\end{split}\end{align}
Now \eqref{g3.1} and \eqref{g3.2} immediately imply
\begin{align}\begin{split}\label{g3.3a}
\norm{A_1(t) (\vek{u},\vek{w})}_{(\mathscr{V}_p)^\ast} \leq c_{\gre,p,\abs{\Omega}}\big(1+\norm{(\vek{u},\vek{w})}_{\mathscr{V}_p}^{p-1} \big)
\end{split}\end{align}
From \eqref{Nbeschr} and H{\"o}lder's inequality we are able to derive that
\begin{align}\begin{split}\label{g3.3b}
\norm{A_3(t) (\vek{u},\vek{w})}_{(\mathscr{V}_p)^\ast} &\leq \sup\limits_{\norm{(\mu,\nu)}_{\mathscr{V}_p}\leq1  } c_{\gre} \big(1+\norm{\vek{w}}_{W_0^{1,p}}^{p-1} \big) \norm{\tpsi}_{W_0^{1,p}} \\
&\leq c_{\gre} \big(1+\norm{(\vek{u},\vek{w})}_{\mathscr{V}_p}^{p-1} \big).
\end{split}\end{align}
To treat the convective terms, we argue the same way as in \cite{br-hirano}. Using H{\"o}lder's inequality 
we immediately get
\begin{align*}\begin{split}
    \norm{A_2(t) (\vek{u},\vek{w})}_{(\mathscr{V}_p)^\ast} &\leq c\,
    \norm{\vek{u}}^2_{L^{2p'}},
    \\
    \norm{A_4(t) (\vek{u},\vek{w})}_{(\mathscr{V}_p)^\ast} &\leq c\,
    \norm{\vek{u}}_{L^{2p'}}\norm{\vek{w}}_{L^{2p'}}\leq c\,\max
    \Big\lbrace \norm{\vek{u}}^2_{L^{2p'}},
      \norm{\vek{w}}^2_{L^{2p'}} \Big\rbrace.
\end{split}\end{align*}
If $p \in [\frac{11}{5},3)$ we use for $r=\frac{12p}{-5p^2+17p-6}$ the H{\"o}lder interpolation
\begin{align}\label{g3.3} \norm{\vek{u}}_{L^r} \leq  \norm{\vek{u}}_{L^2}^{\frac{3-p}{2}} \norm{\vek{u}}_{L^\frac{3p}{3-p}}^\frac{p-1}{2}. \end{align}
Since for these $p $ there holds $2p'\leq r$, we conclude that
\begin{align}\begin{split}\label{g3.4}
\norm{A_2(t) (\vek{u},\vek{w})}_{(\mathscr{V}_p)^\ast} \leq c\,
\norm{\vek{u}}^{3-p}_{L^2}\norm{\vek{u}}_{V_p}^{p-1}\leq c\,
\norm{(\vek{u},\vek{w})}^{3-p}_{\mathscr{H}}\norm{(\vek{u},\vek{w})}_{\mathscr{V}_p}^{p-1}. 
\end{split}\end{align}
Moreover, \eqref{g3.3} also holds for $\vek{w}$, so we also get
\begin{align}\begin{split}\label{g3.5}
 \norm{A_4(t) (\vek{u},\vek{w})}_{(\mathscr{V}_p)^\ast} \leq c\,
 \norm{(\vek{u},\vek{w})}^{3-p}_{\mathscr{H}}\norm{(\vek{u},\vek{w})}_{\mathscr{V}_p}^{p-1}. 
\end{split}\end{align}
If, on the other hand, $p\geq 3$, we get $p'\leq\frac{3}{2}$. This implies that $2<2p'\leq 3$ always holds and therefore it is sufficient to use the interpolation
\begin{align}
\norm{\vek{u}}_{L^3} \leq  \norm{\vek{u}}_{L^2}^{\frac{3-p}{2}} \norm{\vek{u}}_{L^\frac{6(p-1)}{3p-5}}^\frac{p-1}{2}
\end{align}
to get \eqref{g3.4} and \eqref{g3.5} also for $p\geq 3$.
Now \eqref{g3.3a}, \eqref{g3.3b}, \eqref{g3.4} and \eqref{g3.5} imply
\begin{align*}
  \begin{split}
    \norm{A(t) (\vek{u},\vek{w})}_{(\mathscr{V}_p)^\ast} \leq
    C_{\gre,p,|\Omega|}\big (1+\norm{(\vek{u},\vek{w})}_{\mathscr{V}_p}^{p-1}
    \big) + c\, \norm{(\vek{u},\vek{w})}^{3-p}_{\mathscr{H}}\norm{(\vek{u},\vek{w})}_{\mathscr{V}_p}^{p-1}.
  \end{split}
\end{align*}
so that all the assumptions of Theorem \ref{ResultatHirano} are satisfied.
\end{Bew}

\section{Approximative solutions of the system for $\frac{6}{5}<p\leq\frac{11}{5}$}\label{secapprox}

In this section we prove the solvability of an appropriate
approximation of our system \eqref{system}.  The approximate system
arises by adding two terms which are monotone but provide a
better coercivity than the terms induced by $\grs$ and $\grn$.  To
solve this problem we again use Theorem \ref{ResultatHirano}.
\begin{Sa}\label{approxloesung}
  Let $p\in (\frac{6}{5},\frac{11}{5}]$, $T \in (0,\infty)$ and
  $\Omega\subset\R^3$ be a bounded domain with
  Lipschitz-boundary. Assume that $\grs$ satisfies Assumption
  \ref{VssS}, that $\grn $ satisfies Assumption \ref{VssN} and that
  $ \gre$ satisfying Assumption \ref{VssE} is given. Let
  ${\vvv_0\in H(\Omega)}$, $\om_0 \in L^2(\Omega)$ and
  $\f,\l \in L^{p'}(\Omega_T)$ be given.  Then for any
  $q \in (\frac{11}{5},3)$, $M\in \Nat$ there exists
  $(\vvv, \om) \in W_q^1(0,T;\vqom,\hhom) \times
  W_q^1(0,T;\wqom,\lzweiom )$
  satisfying for all  $(\tphi, \tpsi) \in L^q(0,T;V_q(\Omega))\times
  L^q(0,T;W_0^{1,q}(\Omega))$
  \begin{align}
    \begin{split}\label{approxgl}
      &\int_0^T \langle \frac{d\vvv}{dt}(t),\tphi(t)
      \rangle_{V_q}+\int_{\omt}\sdre:\nabla\tphi-\int_{\omt}\vvv\otimes\vvv:\nabla\tphi
      \\ 
      &\quad+\int_0^T \langle \frac{d\om}{dt}(t),\tpsi(t)
      \rangle_{W_0^{1,q}} +\int_{\omt}
      \grn(\nabla\om,\gre):\nabla\tpsi-
      \int_{\omt}\om\otimes\vvv:\nabla\tpsi
      \\ 
      &\quad +\frac{1}{M}\int_{\Omega_T} \abs{\dv}^{q-2} \dv:\grd\tphi
      + \frac{1}{M}\int_{\Omega_T} \abs{\nabla\om}^{q-2}
      \nabla\om:\nabla\tpsi 
      \\ 
      &=\int_{\omt} \f\cdot\tphi +\int_{\omt} \l \cdot
      \tpsi-\int_{\omt} \big (\levy:\sdre \big)\cdot \tpsi
    \end{split}
  \end{align}
  as well as $\vvv(0)=\vvv_0$ and $\om(0)=\om_0$.
\end{Sa}
\begin{Bew}
To apply Theorem \ref{ResultatHirano} we again have to work on a product space. 
Let $\mathscr{V}_q:= V_q(\Omega) \times W_0^{1,q}(\Omega)$,
$\mathscr{H}:= H(\Omega) \times L^2(\Omega)$ and $\mathscr{Z}:= V^3(\Omega)\times W_0^{3,2}(\Omega)$.
Since $q\in (\frac{11}{5},3)$ we get similarly to Section \ref{seceasy} that 
$\mathscr{V}_q,\mathscr{H},\mathscr{Z}$ satisfy the assumptions on the function spaces, which 
are required in Theorem \ref{ResultatHirano}. The operators 
$A(t),A_i(t):\mathscr{V}_q \rightarrow (\mathscr{V}_q)^\ast$, $i=1,\ldots,6$, are defined in \eqref{defOperatoren}$_{1-4}$ and by
\begin{align*}\begin{split}
\left\langle A_5(t)(\vek{u},\vek{w}),(\tphi,\tpsi) \right\rangle&:=\frac{1}{M}\int_{\Omega_T} \abs{\grd\vek{u}}^{q-2} \grd\vek{u}:\grd\tphi,      \\
\left\langle A_6(t)(\vek{u},\vek{w}),(\tphi,\tpsi) \right\rangle&:=\frac{1}{M}\int_{\Omega_T} \abs{\nabla\vek{w}}^{q-2} \nabla\vek{w}:\nabla\tpsi,       \\
A(t)&:= \sum\limits_{i=1}^6 A_i(t),
\end{split}\end{align*}
where $\langle \cdot,\cdot\rangle$ now denotes in all cases the duality pairing between $\mathscr{V}_q$ and $(\mathscr{V}_q)^\ast$.
Again we have to verify that (A1)--(A4) in Theorem
\ref{ResultatHirano} are satisfied.

(A1) From the theory of Nemyckii operators as well as \eqref{Sbeschr}, \eqref{Nbeschr} and 
$E(t)\in L^\infty(\Omega)$ we immediately derive that $A_1(t), A_3(t), A_5(t)$ and $A_6(t)$ 
are continuous operators.
The monotonicity of $A_1(t)$ and $A_3(t)$ is again provided by \eqref{Smonoton} and \eqref{Nmonoton},
whereas $A_5(t)$ and $A_6(t)$ are classical examples of monotone operators (cf.~\cite{lions-quel}, \cite{GGZ}).
Since $A_2(t)$ and $A_4(t)$ are again strongly continuous operators, we see that $A(t)$ is 
pseudomonotone for almost every $t\in [0,T]$.

(A2) This follows in the same way as in the proof of Theorem \ref{easycase}.

(A3) From \eqref{g3.0a} and \eqref{g3.0} follows
\begin{align*}\begin{split}
\left\langle A_2(t)(\vek{u},\vek{w}),\vek{u},\vek{w}) \right\rangle =
\left\langle A_4(t)(\vek{u},\vek{w}),\vek{u},\vek{w}) \right\rangle&=0 ,\\
\left\langle A_1(t)(\vek{u},\vek{w}),(\vek{u},\vek{w}) \right\rangle+\left\langle A_3(t)(\vek{u},\vek{w}),(\vek{u},\vek{w}) \right\rangle &\geq -c(\gre,\Omega).
\end{split}\end{align*}
The definitions of $A_5(t)$ and $A_6(t)$ immediately imply
\begin{align*}\begin{split}
\left\langle A_5(t)(\vek{u},\vek{w}),(\vek{u},\vek{w}) \right\rangle_{\mathscr{V}_q}+\left\langle A_6(t)(\vek{u},\vek{w}),(\vek{u},\vek{w}) \right\rangle_{\mathscr{V}_q}
&= \frac{1}{M} \norm{\vek{u}}_{V_q}^q + \frac{1}{M} \norm{\vek{w}}_{W_0^{1,q}}^q \\
&\geq \frac{c}{M} \norm{(\vek{u},\vek{w})}_{\mathscr{V}_q}^q,
\end{split}\end{align*}
so that (A3) is satisfied.

(A4) To treat $A_2(t)$ and $A_4(t)$ we can use the same argumentation as in Section \ref{seceasy}, if we replace $p$ by $q$. Therefore we get
\begin{align}\begin{split}\label{g4.1}
\norm{A_2(t) (\vek{u},\vek{w})}_{(\mathscr{V}_q)^\ast} +\norm{A_4(t) (\vek{u},\vek{w})}_{(\mathscr{V}_q)^\ast} \leq c\, \norm{(\vek{u},\vek{w})}^{3-q}_{\mathscr{H}}\norm{(\vek{u},\vek{w})}_{\mathscr{V}_q}^{q-1}.
\end{split}\end{align}
Using \eqref{g3.3a}, \eqref{g3.3b}, the continuous embedding $W_0^{1,q} \hookrightarrow W_0^{1,p}$ and Young's inequality, we get
\begin{align}\begin{split}\label{g4.2}
\norm{A_1(t) (\vek{u},\vek{w})}_{(\mathscr{V}_q)^\ast}
&\leq c_{\gre,p,q,\abs{\Omega}}\big (1+\norm{(\vek{u},\vek{w})}_{\mathscr{V}_q}^{q-1} \big),\\
\norm{A_3(t) (\vek{u},\vek{w})}_{(\mathscr{V}_q)^\ast}  &\leq
c_{\gre,\abs{\Omega}} \big (1+\norm{(\vek{u},\vek{w})}_{\mathscr{V}_q}^{q-1} \big).
\end{split}\end{align}
H{\"o}lder's inequality yields
\begin{align}\begin{split}\label{g4.4}
\norm{A_5(t) (\vek{u},\vek{w})}_{(\mathscr{V}_q)^\ast} +\norm{A_6(t) (\vek{u},\vek{w})}_{(\mathscr{V}_q)^\ast}
&\leq \frac{1}{M} \norm{\grd\vek{u}}_{L^q}^{q-1} +\frac{1}{M} \norm{\nabla\vek{w}}_{L^q}^{q-1} \\
&\leq \frac{c}{M} \norm{(\vek{u},\vek{w})}_{\mathscr{V}_q}^{q-1}.
\end{split}\end{align}
Altogether \eqref{g4.1}, \eqref{g4.2}, \eqref{g4.4} imply that (A4) also holds and Theorem \ref{ResultatHirano} reveals the existence of a solution of the problem in Theorem \ref{approxloesung}. 
\end{Bew}

\section{Proof of the main theorem}\label{secmain}
\begin{Sa}[Main Theorem]\label{maintheo}
  Let $p\in (\frac{6}{5},\frac{11}{5}]$, $T \in (0,\infty)$ and
  $\Omega\subset\R^3$ be a bounded domain with
  Lipschitz-boundary. Assume that $\grs$ satisfies Assumption
  \ref{VssS}, that $\grn $ satisfies Assumption \ref{VssN} and that $ \gre$
  satisfying Assumption \ref{VssE} is given.  Then there exists for all
  $\vvv_0\in H(\Omega)$, $\om_0 \in L^2(\Omega)$ and
  $\f,\l \in L^{p'}(\Omega_T)$ a weak solution
  $(\vvv, \om) \in \big (L^p(0,T; V_p(\Omega))\cap L^\infty(0,T;H(\Omega)) \big) \times
  \big (L^p(0,T; W_0^{1,p}(\Omega))\cap L^\infty(0,T;L^2(\Omega))\big)$
  of the problem \eqref{system}--\eqref{ic} satisfying for all
  $(\tphi, \tpsi) \in
  C^{\infty}_{0,\operatorname{div}}([0,T)\times\Omega) \times  C_0^{\infty}([0,T)\times\Omega) $
  \begin{align}
    \begin{split}\label{maineq}
      &-\int_{\omt} \vvv\cdot \partial_t
      \tphi+\int_{\omt}\sdre:\nabla\tphi-\int_{\omt}\vvv\otimes\vvv:\nabla\tphi
      \\ 
      &\quad-\int_{\omt} \om\cdot \partial_t \tpsi +\int_{\omt}
      \grn(\nabla\om,\gre):\nabla\tpsi-
      \int_{\omt}\om\otimes\vvv:\nabla\tpsi
      \\ 
      &\quad+\int_{\omt} \big (\levy:\sdre\big )\cdot \tpsi  
      \\
      &=\int_{\omt} \f\cdot\tphi +\int_{\omt} \l \cdot
      \tpsi+\int_\Omega \vvv_0\cdot\tphi(0)+\int_\Omega
      \om_0\cdot\tpsi(0)\,.
    \end{split}
  \end{align}
\end{Sa}

\begin{Bew}
For any fixed $q \in (\frac{11}{5},3)$  Theorem \ref{approxloesung} yields that for any $M \in \Nat$ there exist solutions
$(\vvv^M,\om^M ) \in W_q^1(0,T;V_q(\Omega),H(\Omega))\times W_q^1(0,T;W_0^{1,q}(\Omega),L^2(\Omega))$ 
with $\vvv^M(0)=\vvv_0$ and $\om(0)=\om_0$ which solve \eqref{approxgl}. Now for any $t\in [0,T]$ 
we are allowed to test \eqref{approxgl} with $\tphi=\vvv^M\chi_{[0,t]} $ and $\tpsi=\om^M\chi_{[0,t]} $. 
Due to $\operatorname{div} \vvv^M=0$ the convective terms vanish and by using the integration by
parts formula \eqref{partielleInteg} together with \eqref{rulelevy} we obtain
\begin{align}\begin{split}\label{g4.1.1}
&  \frac{1}{2}\Big(\norm{\vvv^M(t)}_{H(\Omega)}^2 - \norm{\vvv_0}_{H(\Omega)}^2 +\norm{\om^M(t)}_{L^2(\Omega)}^2-  \norm{\om_0}_{L^2(\Omega)}^2\Big)\\
&\quad  +\frac{1}{M}\int_0^t \int_{\Omega}\abs{\grd(\vvv^M)}^q+\frac{1}{M}\int_0^t \int_{\Omega}\abs{\nabla \om^M}^q+\int_0^t \int_{\Omega}\grn(\nabla\om^M,\gre):\nabla\om^M\\
&\quad  +\int_0^t \int_{\Omega} \grs\big(\grd\vvv^M,\grr(\vvv^M,\om^M),\gre\big):(\grd\vvv^M+\grr(\vvv^M,\om^M))\\
&= \int_0^t \int_{\Omega} \f\cdot\vv^M+\int_0^t \int_{\Omega}\l\cdot\om^M. \raisetag{36mm}
\end{split}\end{align}
We use the coercicity of $\grs$ in \eqref{Skoerziv} and $\grn$ in \eqref{Nkoerziv}, treat the right-hand
side of \eqref{g4.1.1} with H\"{o}lder's, Poincar\'{e}'s, Korn's and Young's inequality and absorb the
resulting terms with $\vvv^M$ and $\om^M$ in the left-hand side of
\eqref{g4.1.1} to get the a~priori
estimate
\begin{align}\begin{split}\label{apriorimain}
&\norm{\vvv^M}_{L^{\infty}(0,T;H(\Omega))}^2+\norm{\om^M}_{L^{\infty}(0,T;L^2(\Omega))}^2
 +\norm{\vvv^M}_{L^p(0,T;V_p(\Omega))}^p\\
&\quad+\norm{\om^M}_{L^p(0,T;W_0^{1,p}(\Omega))}^p+\frac{1}{M}\norm{\vvv^M}_{L^q(0,T;V_q(\Omega))}^q+\frac{1}{M}\norm{\om^M}_{L^q(0,T;W_0^{1,q}(\Omega))}^q\\
&\leq c(\f,\l,\bE,\Omega,\vvv_0,\om_0). \raisetag{18mm}
\end{split}\end{align}
The growth condition of $\grs$ and $\grn$
 (\eqref{Sbeschr} and \eqref{Nbeschr}) together with the theory of Nemyckii operators and 
 \eqref{apriorimain} yield
\begin{align}\begin{split}\label{SNbound}
\norm{\grs\big(\grd\vv^M,\grr(\vv^M,\om^M),\gre\big)}_{L^{p'}(\omt)}+\norm{\grn(\nabla\om^M,\gre)}_{L^{p'}(\omt)}\leq c.
\end{split}\end{align}
The parabolic interpolation $L^p(0,T;W_0^{1,q}(\Omega))\cap L^\infty (0,T;L^2(\Omega))\hookrightarrow
L^{\frac{5}{3}p}(\Omega_T)$ implies
\begin{align}\label{wirbelbound}
\norm{\vvv^m\otimes \vvv^M}_{L^{\frac{5}{6}p}(\Omega_T)}+\norm{\om^M\otimes\vvv^M}_{L^{\frac{5}{6}p}(\Omega_T)}\leq c.
\end{align}
Furthermore there holds 
\begin{align}\begin{split}\label{hilfstermekonv}
\frac{1}{M} \abs{\grd\vvv^M}^{q-2} \grd\vvv^M \overset{M\rightarrow \infty}{\longrightarrow} &\bfzero \qquad \text{in} \ L^{q'}(\omt),\\
\frac{1}{M} \abs{\nabla\om^M}^{q-2} \nabla\om^M \overset{M\rightarrow \infty}{\longrightarrow} &\bfzero \qquad \text{in} \ L^{q'}(\omt).
\end{split}\end{align}
Now we choose $\sigma \in \R$ which satisfies
\begin{align}\label{wahlsigma}
1<\sigma < \min\{\frac{5p}{6},q'\} \ \text{and} \ 2\sigma>p.
\end{align}
This is always possible,
since the second inequality in \eqref{wahlsigma} only provides a restriction if $p>2$ and in this case we have $\min\{\frac{5p}{6},q'\}>\frac{3}{2}$ so $\sigma=\frac{3}{2}$ satisfies
\eqref{wahlsigma}. The next step in our proof is to show the boundedness of $(\vvv^M)$ in 
$W^{1,p,\sigma}(0,T;V_p(\Omega),V_{\sigma'}(\Omega)^\ast)$ and of $(\om^M)$ in 
$W^{1,p,\sigma}(0,T;W_0^{1,p}(\Omega),W_0^{1,\sigma'}(\Omega)^\ast)$. Let
us firstly treat $(\vvv^M)$. 
First due to $\sigma<\frac{5p}{6}$ and $p\leq \frac{11}{5}$ we have $\sigma'>p$. This together with
$p>\frac{6}{5}$ ensures 
\begin{align*}
V_p(\Omega)\hookrightarrow H(\Omega) \hookrightarrow V_p(\Omega)^\ast \hookrightarrow  V_{\sigma'}(\Omega)^\ast
\end{align*}
with continuous and dense embeddings. Therefore
 $W^{1,p,\sigma}(0,T;V_p(\Omega),V_{\sigma'}(\Omega)^\ast)$ is a Bochner-Sobolev space as introduced 
 in Section \ref{subsecnotation}.
From \eqref{apriorimain} we already know that $\norm{\vvv^M}_{L^p(0,T;V_p(\Omega))}$ is bounded 
 and because of $\sigma<q'$ we have that 
$L^{q'}(0,T; V_q(\Omega)^\ast ) \hookrightarrow L^{\sigma}(0,T; V_{\sigma'}(\Omega)^\ast )$, which allows us
 to interpret $\frac{d\vvv^M}{dt}$ as an element of $L^{\sigma}(0,T; V_{\sigma'}(\Omega)^\ast )$. To show the 
boundedness we use arbitrary ${\tphi \in L^{\sigma'}(0,T; V_{\sigma'}(\Omega))}$ and $\tpsi=\bfzero$ in
 \eqref{approxgl}. With the help of H\"{o}lder's inequality as well as \eqref{SNbound}--\eqref{wahlsigma} we estimate
\begin{align}\label{dvmbound}
\Big\Vert\frac{d\vvv^M}{dt}\Big\Vert_{L^{\sigma}(0,T; V_{\sigma'}(\Omega)^\ast )} \leq c.
\end{align}
For $\om^M$ we proceed analogously. The choice of $\sigma$ ensures that 
\begin{align*}
W_0^{1,p}(\Omega)\hookrightarrow L^2(\Omega) \hookrightarrow W_0^{1,p}(\Omega)^\ast \hookrightarrow W_0^{1,\sigma'}(\Omega)^\ast 
\end{align*} 
with continuous and dense embeddings. In \eqref{apriorimain} we already proved the boundedness of 
$\om^M$ in $L^p(0,T;W_0^{1,p}(\Omega))$ and to get an estimate for the time derivative we test \eqref{approxgl} with $\tphi=\bfzero$
 and arbitrary $\tpsi \in L^{\sigma'}(0,T;W_0^{1,\sigma'}(\Omega))$. Then again H\"{o}lder's inequality, 
 the choice of $\sigma$ and \eqref{SNbound}--\eqref{wahlsigma} imply 
\begin{align}\label{dommbound}
\Big\Vert\frac{d\om^M}{dt}\Big\Vert_{L^{\sigma}(0,T; W_0^{1,\sigma'}(\Omega)^\ast )} \leq c.
\end{align}
From \eqref{apriorimain}, \eqref{dvmbound}, \eqref{dommbound} we get
\begin{align}\label{ommbound}
\norm{\vvv^M}_{W^{1,p,\sigma}(0,T;V_p(\Omega),V_{\sigma'}(\Omega)^\ast)}+\norm{\om^M}_{W^{1,p,\sigma}(0,T;W_0^{1,p}(\Omega),W_0^{1,\sigma'}(\Omega)^\ast)} \leq c.
\end{align}
Now \eqref{apriorimain}, \eqref{SNbound}, \eqref{ommbound}, the Aubin-Lions lemma and 
parabolic interpolation lead to the following convergence results after choosing appropriate subsequences:
\begin{alignat}{2}\raisetag{3cm}\label{konvmain}\begin{aligned}
\vvv^M &\overset{\ast}{\rightharpoonup} \vvv& \quad &\text{in} \ L^{\infty}(0,T;H(\Omega)),\\
\om^M &\overset{\ast}{\rightharpoonup} \om&&\text{in} \ L^{\infty}(0,T;L^2(\Omega)),\\
\vvv^M &\rightharpoonup \vvv&&\text{in} \ W^{1,p,\sigma}(0,T;V_p(\Omega), V_{\sigma'}(\Omega)^{\ast}),\\
\om^M &\rightharpoonup \om&&\text{in} \ W^{1,p,\sigma}(0,T;W_0^{1,p}(\Omega), W_0^{1,\sigma'}(\Omega)^{\ast}),\\
\vvv^M &\rightarrow \vvv&&\text{in} \ L^{2\sigma}(\omt),\\
\om^M &\rightarrow \om&&\text{in} \ L^{2\sigma}(\omt),\\
\vvv^M\otimes\vv^M &\rightarrow \vvv\otimes\vv&&\text{in} \ L^{\sigma}(\omt),\\
\om^M\otimes\vvv^M &\rightarrow \om\otimes\vvv&&\text{in} \ L^{\sigma}(\omt),\\
\grs\big(\grd\vvv^M,\grr(\vvv^M,\om^M),\gre\big) &\rightharpoonup \widehat{\grs} &&\text{in} \ L^{p'}(\omt),\\
\grn(\nabla\om^M,\gre) &\rightharpoonup \widehat{\grn} &&\text{in} \ L^{p'}(\omt).
\end{aligned}
\end{alignat}
Here we also made use of our choice of $\sigma$ since the Aubin-Lions lemma and 
parabolic interpolation imply $\vvv^M \rightarrow \vvv,\ \om^M \rightarrow \om$ in $L^s(\omt)$ for any 
\mbox{$1\leq s < \frac{5}{3}p$.}
In particular $\vvv$ and $\om$ belong to the required function spaces
in Theorem \ref{maintheo}.

To derive our limit equation we test \eqref{approxgl}, which is solved for any $M\in\Nat$ by $\vvv^M$ and
 $\om^M$, with arbitrary $\tphi \in C^{\infty}_{0,\operatorname{div}}([0,T)\times\Omega)$ and $\tpsi \in
  C_0^{\infty}([0,T)\times\Omega)$ and use the integration by parts formula. We get
\begin{align}\begin{split}\label{g4.1.2}
&-\int_{\omt}\vvv^M\cdot \partial_t \tphi+\int_{\omt}\grs\big(\grd\vvv^M,\grr(\vvv^M,\om^M),\gre\big):\nabla\tphi-\int_{\omt}\vvv^M\otimes\vvv^M:\nabla\tphi\\
&\quad-\int_{\omt}\om^M\cdot \partial_t \tpsi +\int_{\omt} \grn(\nabla\om^M,\gre):\nabla\tpsi- \int_{\omt}\om^M\otimes\vvv^M:\nabla\tpsi\\
&\quad+\frac{1}{M}\int_{\omt}\abs{\grd\vvv^M}^{q-2}\grd\vvv^M:\grd\tphi+\frac{1}{M}\int_{\omt}\abs{\nabla\om^M}^{q-2}\nabla\om^M:\nabla\tpsi\\
&\quad+\int_{\omt}\big (\levy:
\grs\big(\grd\vvv^M,\grr(\vvv^M,\om^M),\gre\big)\big )\cdot \tpsi\\
&=\int_{\omt} \f\cdot\tphi +\int_{\omt} \l \cdot \tpsi +\int_\Omega \vvv^M(0)\cdot\tphi(0)+\int_\Omega \om^M(0)\cdot\tpsi(0). \raisetag{20mm}
\end{split}\end{align}
Since for any $M\in\Nat$ there holds $\vvv^M(0)=\vvv_0$ in $H(\Omega)$ and $\om^M(0)=\om_0$ in $L^2(\Omega)$
the convergences in \eqref{hilfstermekonv} and \eqref{konvmain} allow us to pass to the limit in every
 term and we conclude 
\begin{align}\begin{split}\label{grenzwertgl}
&-\int_{\omt}\vvv\cdot \partial_t \tphi+\int_{\omt}\widehat{\grs}:\nabla\tphi-\int_{\omt}\vvv\otimes\vvv:\nabla\tphi\\
&\quad-\int_{\omt}\om\cdot \partial_t \tpsi +\int_{\omt} \widehat{\grn}:\nabla\tpsi- \int_{\omt}\om\otimes\vvv:\nabla\tpsi+\int_{\omt} (\levy:\widehat{\grs})\cdot \tpsi\\
&=\int_{\omt} \f\cdot\tphi +\int_{\omt} \l \cdot \tpsi+\int_\Omega \vvv_0\cdot\tphi(0)+\int_\Omega \om_0\cdot\tpsi(0) 
\end{split}\end{align}
holds for for all $\tphi \in C^{\infty}_{0,\operatorname{div}}([0,T)\times\Omega)$ and all
$\tpsi \in C_0^{\infty}([0,T)\times\Omega)$.

So it remains to proof that $\widehat{\grs}=\sdre$
and $\widehat{\grn}=\nome$
a.e.~in $\omt$
to finish the proof of Theorem \ref{maintheo}. To this end we use the
two Lipschitz truncation results in Section \ref{sec:aux}.

\smallskip 
We start with the proof of $\widehat{\grs}=\sdre$ a.e.~in $\omt$. We
set $\tpsi=\bfzero$ and choose
$\tphi \in C_{0,\operatorname{div}}^\infty(\omt)$ arbitrarily in
\eqref{g4.1.2} and \eqref{grenzwertgl} and subtract these two
equations. Thus we obtain for any $\tphi \in C_{0,\operatorname{div}}^\infty(\omt)$
\begin{align}
  \begin{split}\label{g4.1.3}
    -\int_{\omt}(\vvv^M-\vvv)\cdot \partial_t \tphi= &\int_{\omt}(
    \widehat{\grs}- \grs\big(\grd\vvv^M,\grr(\vvv^M,\om^M),\gre\big)
    ):\nabla\tphi 
    \\
    &+\int_{\omt} (\vvv^m\otimes\vvv^M - \vvv\otimes\vvv- \frac{1}{M}
    \abs{\grd\vvv^M}^{q-2}\grd\vvv^M ):\nabla\tphi .\raisetag{15mm}
  \end{split}
\end{align}
Now we argue by contradiction. Assume that there exists a set $\mathcal{M} \subseteq \omt$, which satisfies
$\abs{\mathcal{M}}\geq 2\delta$ for some $\delta>0$, so that almost everywhere in $\mathcal{M}$ there holds
$\widehat{\grs}\neq \sdre$. For $\epsilon>0$ we define 
\[\Omega_{\epsilon}:=\{x\in \Omega   \fdg \operatorname{dist}(\partial\Omega,x)\geq \epsilon \} 
\ \textnormal{and} \ 
\Omega_{T,\epsilon}:=[\epsilon,T-\epsilon]\times \Omega_\epsilon.\]
Clearly we can choose $\epsilon>0$ sufficiently small so that $\abs{\mathcal{M}\cap \Omega_{T,\epsilon}} \geq\delta$. Since $\Omega_{T,\epsilon}$ is compact there exists $n\in\Nat$ and $(t_i,x_i)_{1\leq i\leq n} \subset\Omega_{T,\epsilon}$ such that
\[
\Omega_{T,\epsilon} \subset \bigcup_{i=1}^n(t_i-\tfrac{\epsilon}{8},t_i+\tfrac{\epsilon}{8})\times B_{\frac{\epsilon}{8}}(x_i).
\]
This in turn implies the existence of $j \in \{1,\ldots,n\}$ so that
\begin{align*}
\abs{ \mathcal{M}\cap \big (
  (t_j-\tfrac{\epsilon}{8},t_j+\tfrac{\epsilon}{8})\times
  B_{\frac{\epsilon}{8}}(x_j)  \big )   }\geq \tfrac{\delta}{n}>0.
\end{align*}
Therefore it is sufficient to prove that  $\widehat{\grs}=\sdre$ holds
a.e.~in $(t_j-\tfrac{\epsilon}{8},t_j+\tfrac{\epsilon}{8})\times
B_{\frac{\epsilon}{8}}(x_j)   $ to achieve a  contradiction. 
To this end we set $I_0:= (t_j -\epsilon,t_j+\epsilon)$, $B_0:= B_{\epsilon}(x_j)$, $Q_0:=I_0\times B_0$
and define
\begin{align*}\begin{split}
\bfu_M&:= (\vvv^M-\vvv)\chi_{Q_0},\\
\bfG_{1,M}&:= \big (\widehat{\grs}- \grs\big(\vek{D}\vvv^M,\grr(\vvv^M,\om^M),\gre\big)\big)\chi_{Q_0},\\
\bfG_{2,M}&:= (\vvv^M\otimes\vvv^M  -   \vvv\otimes\vvv   -    \frac{1}{M}\abs{\grd\vvv^M}^{q-2}\grd\vvv^M                 ) \chi_{Q_0},\\
\bfG_M&:= \bfG_{1,M}+\bfG_{2,M}.
\end{split}\end{align*}
With these definitions we get from \eqref{g4.1.3}, by restricting the
test functions, that for any $\boldsymbol{\xi} \in C_{0,\Div}^{\infty}(Q_0)$ 
\begin{align}\begin{split}\label{eq:aa}
-\int_{Q_0} \mathbf{u}_M\cdot \partial_t \boldsymbol{\xi}  =\int_{Q_0} \mathbf{G}_M: \nabla \boldsymbol{\xi}.
\end{split}\end{align}
The functions $\bfu_M$ and $\bfG_M$ satisfy the assumptions of Theorem
\ref{solpara} because of \eqref{hilfstermekonv} and \eqref{konvmain}.
Thus Theorem \ref{solpara} yields that there exist double-sequences
$(\lambda_{M,k})$ and $(\mathcal O_{M,k})$ such that for every
$k\geq k_0$ there holds $2^{2^k} \leq \lambda_{M,k} \leq 2^{2^{k+1}}$
and $\abs{\mathcal{O}_{M,k}}\leq 2^{-k}$. If we choose
$\zeta \in C_0^{\infty}(\frac{1}{6}Q_0)$ with
$\chi_{\frac{1}{8}Q_0}\leq \zeta \leq \chi_{\frac{1}{6}Q_0} $ as well
as\footnote{Since $\vvv \in L^p(0,T;V_p(\Omega))$ and
  $\om \in L^p(0,T;W_0^{1,p}(\Omega))$ the growth condition of $\grs$
  ensures that $\sdre \in L^{p'}(\omt)$. }
$\bfK:= \sdre-\widehat{\grs}\in L^{p'}(\frac{1}{6}Q_0)$ we get from
Corollary \ref{corsolpara} that
\begin{align}\begin{split}\label{g4.1.4}
\limsup_{M\rightarrow \infty} \Big|  \int_{ \frac{1}{6}Q_0 } (\grs\big(\vek{D}\vvv^M,\grr(\vvv^M,\om^M),\gre\big) -\sdre):&\\
\nabla (\vvv^M-\vvv)\zeta \chi_{\mathcal{O}_{M,k}^\complement}  \Big| &\leq 2^{\frac{-k}{p}}.
\end{split}\end{align}
Unfortunately the integrand has no sign, since we can't use the monotonicity of $\grs$ in \eqref{Smonoton}. 
On the other hand $\grs\big(\vek{D}\vvv^M,\grr(\vvv^M,\om^M),\gre\big) -\sdre$ is bounded in $L^{p'}(\Omega_T)$ and 
due to the choice of $\sigma$ in \eqref{wahlsigma}, \eqref{konvmain} and \eqref{levyabsch} we can deduce that
 $\levy\cdot\om^M \rightarrow \levy \cdot\om$ in $L^p(\Omega_T)$. Therefore 
\begin{align}\begin{split}\label{g4.1.5}
\limsup_{M\rightarrow \infty} \Big|  \int_{ \frac{1}{6}Q_0 } (\grs\big(\vek{D}\vvv^M,\grr(\vvv^M,\om^M),\gre\big) -\sdre):&\\
(\levy\cdot\om^M - \levy \cdot\om)\zeta \chi_{\mathcal{O}_{M,k}^\complement}  \Big| &= 0.
\end{split}\end{align}
Form \eqref{g4.1.4}, \eqref{g4.1.5} and the definition of $\grr$ we conclude that
\begin{align}\begin{split}\label{g4.1.6}
\limsup_{M\rightarrow \infty} \Big|  \int_{ \frac{1}{6}Q_0 }
\big (\grs\big(\vek{D}\vvv^M,\grr(\vvv^M,\om^M),\gre\big) -\sdre\big ):&\\
\big (\grd\vvv^M + \grr(\vvv^M,\om^M) - \grd\vvv -\grr(\vvv,\om)     \big )\zeta \,\chi_{\mathcal{O}_{M,k}^\complement}  \Big| &\leq c\, 2^{\frac{-k}{p}}.\raisetag{15mm}
\end{split}\end{align}
Now this term has a sign due to \eqref{Smonoton}. From now on we use the abbreviation 
\begin{align*}
  \mathcal{S}^M:= \big
  (\grs\big(\vek{D}\vvv^M,\grr(\vvv^M,\om^M),\gre\big) -\sdre\big ): 
  \\ 
  \big(\grd\vvv^M + \grr(\vvv^M,\om^M) - \grd\vvv -\grr(\vvv,\om)\big).
\end{align*}
Note that from \eqref{konvmain} we are able to conclude that $\mathcal{S}^M$ is bounded in $L^1(\omt)$.
Since $\mathcal{S}^M\geq 0$ the expression
$(\mathcal{S}^M)^{\frac{1}{2}}$ is well defined. From H\"{o}lder's inequality, $\abs{\zeta}\leq 1$ and $\abs{\mathcal{O}_{M,k}}\leq c\, 2^{-k}$ we easily get
\begin{align*}\begin{split}
\limsup_{M\rightarrow \infty}\Big|  \int_{ \frac{1}{6}Q_0 }
(\mathcal{S}^M)^{\frac{1}{2}} \zeta \chi_{\mathcal{O}_{M,k}} \Big|
\leq \limsup_{M\rightarrow \infty} \norm{ \mathcal{S}^M}_{L^1(
  \frac{1}{6}Q_0 )}^{\frac{1}{2}} \abs{  \mathcal{O}_{M,k}
}^{\frac{1}{2}} \leq c\, 2^{\frac{-k}2}, 
\end{split}\end{align*}
which together with \eqref{g4.1.6} implies
\begin{align}\begin{split}\label{g4.1.8}
\limsup_{M\rightarrow \infty} \int_{ \frac{1}{6}Q_0 }
\min\left\lbrace\mathcal{S}^M ,( \mathcal{S}^M)^{\frac{1}{2}}
\right\rbrace \zeta \leq c\,\min \Big
\{2^{\frac{-k}{p}},2^{\frac{-k}{2}}\Big \}.
\end{split}\end{align}
Since this holds for any $k\geq k_0$ and since $\zeta=1$ on $\frac{1}{8}Q_0$ we get 
\begin{align}\label{g4.1.9}
\mathcal{S}^M \rightarrow 0 \qquad \text{almost everywhere in} \ \tfrac{1}{8}Q_0,
\end{align}
at least for a not relabeled subsequence. Now we define $\vek{T}:\R^{3\times3} \rightarrow \R^{3\times3}$ by
$\vek{T}(\bfA)=\grs\big(\bfA^{\operatorname{sym}},\bfA^{\operatorname{skew}},\gre \big)$. Due to our Assumptions \ref{VssS} and \ref{VssE} on $\grs$ and $\gre$ it is clear that $\vek{T}$ is continuous and strictly monotone.
So if we set $\boldsymbol{\eta}^M(t,x):= \grd\vv^M(t,x)+\grr(\vv^M,\om^M)(t,x) $ and $\boldsymbol{\eta}(t,x):= \grd\vv(t,x)+\grr(\vv,\om)(t,x) $ the equation \eqref{g4.1.9} reads 
\begin{align*}\begin{split}
\big(\vek{T}(\boldsymbol{\eta}^M(t,x) )-\vek{T}(\boldsymbol{\eta}(t,x) )\big):\big(\boldsymbol{\eta}^M(t,x) -\boldsymbol{\eta}(t,x) \big) \rightarrow 0
\end{split}\end{align*}
almost everywhere in $\frac{1}{8}Q_0$. Thus \cite[Lemma 6]{DMM} implies 
\begin{align*}\begin{split}
\boldsymbol{\eta}^M(t,x) \rightarrow \boldsymbol{\eta}(t,x)
\end{split}\end{align*}
a.e.~in $\frac{1}{8}Q_0$. Since we already know from \eqref{konvmain} that $\levy\cdot\om^M
\rightarrow \levy\cdot \om$ holds a.e.~in $\frac{1}{8}Q_0$, we conclude that 
$\grd\vvv^M \rightarrow \grd\vvv$ and
$\grw\vvv^M \rightarrow \grw\vvv$ also holds a.e.~in $\frac{1}{8}Q_0$. Since $\grs$ 
is continuous this implies
\begin{align}\begin{split}\label{konvae}
\grs\big(\grd\vvv^M,\grr(\vvv^M,\om^M),\gre     \big) \rightarrow \sdre \quad \text{a.e. in } \tfrac{1}{8}Q_0.
\end{split}\end{align}
Since weak and a.e.~limits coincide (cf.~\cite{Hew-strom-65}) we conclude from \eqref{konvmain} and \eqref{konvae} that $\sdre= \widehat{\grs}$ a.e.~in 
$\tfrac{1}{8}Q_0 =(t_j-\tfrac{\epsilon}{8},t_j+
   \tfrac{\epsilon}{8})\times B_{\frac{\epsilon}{8}}(x_j), $ 
which gives the desired contradiction and proofs 
\begin{align*}\begin{split}
\sdre= \widehat{\grs} \qquad \text{almost everywhere in } \omt.
\end{split}\end{align*}
\newline
So the remaining step in the proof of the main Theorem \ref{maintheo} is to prove that $\widehat{\grn}=\nome$ 
holds a.e.~in $\omt$. We start by subtracting the equations \eqref{g4.1.2} and \eqref{grenzwertgl} 
one from another. If we set $\tphi=\bfzero$ we get for any $\tpsi \in
C_0^\infty(\omt)$ 
\begin{align}\begin{split}\label{g4.1.10}
-\int_{\omt} (\om^M-\om)\cdot \partial_t\tpsi=& \int_{\omt} (\widehat{\grn}-\grn(\nabla \om^M,\gre)):\nabla\tpsi\\
&+\int_{\omt} ( \om^m\otimes \vvv^M -\om\otimes\vvv -\frac{1}{M} \abs{\nabla \om^M}^{q-2}\nabla \om^M ):\nabla\tpsi\\
&+\int_{\omt} \Big( \levy:\big(
\grs\big(\vek{D}\vvv^M,\grr(\vvv^M,\om^M),\gre\big) -\widehat{\grs}
\big)   \Big)\cdot \tpsi .\raisetag{27mm}
\end{split}\end{align}
The ideas we want to use have been developed in \cite{die-ru-wolf} but
there are some differences. The first difference is that
\eqref{g4.1.10} holds for any $\tpsi \in C_0^\infty(\omt)$ not just
for $\tpsi \in C_{0,\operatorname{div}}^\infty(\omt)$ as in
\cite{die-ru-wolf}. Therefore we don't need to construct a local
pressure which had to be done in \cite{die-ru-wolf}. The second
difference is that in our equation \eqref{g4.1.10} a term of lower
order, namely
$\int_{\omt} \big( \levy:\big(
\grs\big(\vek{D}\vvv^M,\grr(\vvv^M,\om^M),\gre\big) -\widehat{\grs}
\big) \big)\cdot \tpsi$,
appears. This is why we had to prove a slightly generalized Lipschitz
Truncation in Theorem \ref{paralip}.  We define
\begin{align}\begin{split}\label{g4.1.11}
\uuum&:=\om^M-\om,\\
\vek{H}_{1,M}&:= \widehat{\grn}-\grn(\nabla \om^M,\gre),\\
\vek{H}_{2,M}&:= \om^M\otimes \vvv^M -\om\otimes\vvv -\frac{1}{M} \abs{\nabla \om^M}^{q-2}\nabla \om^M,\\
\vek{H}_{M}&:=\vek{H}_{1,M}+\vek{H}_{2,M},\\
\vek{k}_M&:= \levy:\big( \grs\big(\vek{D}\vvv^M,\grr(\vvv^M,\om^M),\gre\big) -\widehat{\grs}    \big), 
\end{split}\end{align}
so that \eqref{g4.1.10} reads  
\begin{align}\begin{split}\label{g4.1.12}
-\int_{\omt} \uuum\cdot\partial_t \tpsi = \int_{\omt}\vek{H}_{M}:\nabla \tpsi +  \int_{\omt}\vek{k}_M\cdot \tpsi
\end{split}\end{align}
for any $\tpsi \in C_0^\infty(\omt)$ just as in \eqref{g2.3.3} in
Theorem \ref{paralip}. Using a density argument, we conclude that for
any $\tpsi \in L^{\sigma'}(0,T;W_0^{1,\sigma'}(\Omega))$ 
\begin{align}\begin{split}\label{g4.1.14}
\int_0^T \big \langle \frac{d\uuum}{dt}(\tau),\tpsi(\tau)
\big \rangle_{W_0^{1,\sigma'}(\Omega)}= \int_{\omt}\vek{H}_{M}:\nabla
\tpsi +  \int_{\omt}\vek{k}_M\cdot \tpsi .
\end{split}\end{align}
We also already know
from \eqref{hilfstermekonv} and \eqref{konvmain} that
\begin{alignat}{2}\label{g4.1.13}\begin{aligned}
\uuum &\rightharpoonup \bfzero   &&\text{in } W^{1,p,\sigma}(0,T;W_0^{1,p}(\Omega), W_0^{1,\sigma'}(\Omega)^{\ast}),\\
\uuum&\rightarrow \bfzero   \qquad &&\text{in } L^{2\sigma}(\omt),\\
\vek{H}_{1,M}&\rightharpoonup \bfzero   \qquad &&\text{in } L^{p'}(\omt),\\
\vek{H}_{2,M}&\rightarrow \bfzero &&\text{in } L^\sigma(\omt), \\
\vek{k}_M& \rightharpoonup\bfzero &&\text{in } L^{p'}(\omt).
\end{aligned}
\end{alignat}
Now we have to choose similar to \cite{die-ru-wolf} a double sequence $(\lMk)$ and exceptional 
sets $(\emk)$ for which we want to apply Theorem \ref{paralip}. To
this end we set
\begin{align*}\begin{split}
g_M:= \mathcal{M}^\ast(\abs{\nabla \uuum}) + \mathcal{M}^\ast(\abs{\vek{H}_{1,M}})^{\frac{1}{p-1}}+ \mathcal{M}^\ast(\abs{\vek{k}_M})^{\frac{1}{p-1}}.
\end{split}\end{align*}
Here we extended all appearing functions by zero to the whole space. Due to the 
strong-type estimate of the maximal operator $\mathcal{M}^\ast$ we obtain
\begin{align*}\begin{split}
\norm{g_M}_{L^p} &\leq \norm{\mathcal{M}^\ast(\abs{\nabla \uuum})  }_{L^p}+ \norm{\mathcal{M}^\ast(\abs{\vek{H}_{1,M}})}_{L^{p'}}^{\frac{p'}{p}}+ \norm{\mathcal{M}^\ast(\abs{\vek{k}_M})}_{L^{p'}}^{\frac{p'}{p}}\\
&\leq c(p,p')\, \big(\norm{\nabla \uuum  }_{L^p}+ \norm{\vek{H}_{1,M}}_{L^{p'}}^{\frac{p'}{p}}+ \norm{\vek{k}_M}_{L^{p'}}^{\frac{p'}{p}}\big) \leq c.
\end{split}\end{align*}
Therefore we have for any $k\in \Nat$ 
\begin{align*}\begin{split}
c^p \geq \int\limits_{2^{2^k}}^{2^{2^{k+1}}} p \lambda^{p-1} \abs{ \{ \abs{g_M}>\lambda   \}   } d\lambda
&\geq p \int\limits_{2^{2^k}}^{2^{2^{k+1}}} \lambda^{-1} d\lambda \inf_{2^{2^k} \leq \gamma \leq 2^{2^{k+1}}}
\gamma^p \abs{ \{ \abs{g_M}>\gamma   \}   }\\
&\geq p\,2^k \ln(2) \inf_{2^{2^k} \leq \gamma \leq 2^{2^{k+1}}}
\gamma^p \abs{ \{ \abs{g_M}>\gamma   \}   }
\end{split}\end{align*}
so that we are able to choose 
\begin{align}\label{lmkabsch}
\lMk\in \left[2^{2^k},  2^{2^{k+1}}  \right] 
\end{align}
such that
\begin{align}\label{g4.1.15}
\lMk^p \abs{ \{ \abs{g_M}>\lMk   \}  } \leq c\,2^{-k}
\end{align}
holds for any $k,M \in \Nat$. For any $k,M \in \Nat$ we define $G_{M,k}:= \{ \abs{g_M}>\lMk   \}$ and 
$\aMk:=\lMk^{2-p}$. Now \eqref{g4.1.15} reads 
\begin{align}\label{gmkabsch}
\lMk^p \abs{G_{M,k}} \leq c\,2^{-k}.
\end{align}
Since $\frac{\lMk}{\aMk}=\lMk^{p-1}$, we have
\begin{align}\begin{split}\label{gmkobermenge}
G_{M,k}&= \left\{ \mathcal{M}^\ast(\abs{\nabla \uuum}) + \mathcal{M}^\ast(\abs{\vek{H}_{1,M}})^{\frac{1}{p-1}}+ \mathcal{M}^\ast(\abs{\vek{k}_M})^{\frac{1}{p-1}}>\lMk \right\}\\
 &\supset \left\{ \mathcal{M}^\ast (\abs{\nabla \uuum}) >\lMk  \right\} \cup \left\{ \amk\mathcal{M}^\ast (\abs{\vek{H}_{1,M}}) >\lMk   \right\} \\
&\qquad\cup  \left\{ \aMk  \mathcal{M}^\ast (\abs{\vek{k}_M}) > \lMk  \right\}.
\end{split}\end{align}
Next we define 
\begin{align}\label{fmkgl}
 F_{M,k}:= \left\{ \mathcal{M}^\ast(\abs{\vek{H}_{2,M}})>\lMk^{p-1}\right\} =
  \left\{ \aMk\mathcal{M}^\ast(\abs{\vek{H}_{2,M}})>\lMk\right\}.
 \end{align}  
Then the weak-type estimate of $\mathcal{M}^\ast$ and \eqref{g4.1.13} imply
\begin{align}\label{fmkab}
\abs{F_{M,k}}\leq c (\lMk^{p-1})^{-\sigma} \norm{\vek{H}_{2,M}}_{L^\sigma}^\sigma \stackrel{M\rightarrow \infty}{\longrightarrow}0
\end{align}     
for any fixed $k\in\Nat$. Since $\mathcal{M}^\ast$ is subadditive we conclude from 
\eqref{gmkobermenge} and \eqref{fmkgl} that 
\begin{align}\begin{split}\label{gmkfmkobermenge}
G_{M,k} \cup F_{M,k} \supset  \left\{ \mathcal{M}^\ast(\abs{\nabla \uuum}) + \aMk\mathcal{M}^\ast(\abs{\vek{H}_{M}})+ \aMk\mathcal{M}^\ast(\abs{\vek{k}_M})>4\lMk \right\}.
\end{split}\end{align}
Moreover, from the fact that $\{\mathcal{M}(f)>\beta\} $ is an open set for any $f \in L^s$ and $\beta>0$ it
is easy to prove that $G_{M,k}$ and $F_{M,k}$ are open sets as well. We also define for each 
$M,k \in \Nat$ the set
\begin{align}\label{hmkdef}
H_{M,k}:=\left\lbrace \mathcal{M}^\ast(\abs{\uuum})>1 \right\rbrace
\end{align}
so that the weak-type estimate for $\mathcal{M}^\ast$ and \eqref{g4.1.13} imply
\begin{align}\label{hmkabsch}
\abs{H_{M,k}} \leq c \norm{\uuum}_{L^{2\sigma}}^{2\sigma} \stackrel{M\rightarrow \infty}{\longrightarrow} 0.
\end{align}
Now we can define our exceptional set
\begin{align}\label{emkdef}
\emk:= ( G_{M,k} \cup  F_{M,k} \cup H_{M,k}) \cap \omt.
\end{align}
Clearly from \eqref{gmkabsch}, \eqref{fmkab} and \eqref{hmkabsch} we get for any fixed $k \in \Nat$
\begin{align}\label{emkAbsch}
\limsup_{M\rightarrow \infty} \lMk^p \abs{\emk} \leq c\,2^{-k}.
\end{align}
We choose an arbitrary cut-off function $\zeta \in C_0^\infty(\omt)$ and define 
the compact set $K:=\operatorname{supp}\zeta$. Due to \eqref{gmkfmkobermenge} and 
\eqref{hmkdef} the set $\emk$ satisfies
\begin{align}\label{emkTeilmenge}
\omt \cap (\mathcal{O}_{4\lMk} \cup \mathcal{U}_1  ) \subset \emk \subset \omt
\end{align}
so that we are able to apply Theorem \ref{paralip} with 
$\Lambda=4\lMk$, $\alpha=\aMk$, $\bfu=\uuum$, $\vek{H}=\vek{H_M}$, 
$\vek{k}=\vek{k}_M$, $E=\emk$ and $K= \operatorname{supp}\zeta$ for any 
$M,k \in \Nat$. To ensure a better readability we denote $\trunckm:=\mathcal{T}_{\emk}^{\aMk}$. 
Due to Theorem \ref{paralip} (i) the function $(\trunckm\uuum )\zeta$ is an 
admissible test function for \eqref{g4.1.14} and due to Theorem \ref{paralip} 
(iv) we have
\begin{align*}\begin{split}
 \int_0^T &\big\langle \frac{d\uuum(\tau)}{dt}  , (\trunckm\uuum(\tau))\zeta(\tau)  \big\rangle_{W_0^{1,\sigma}} d\tau\\
 = \frac{1}{2}&\int_{\omt}(\abs{\trunckm\uuum}^2 - 2\uuum\cdot \trunckm\uuum      )\partial_t \zeta\\
 &+\int_{E_{m,k}}(\partial_t \trunckm\uuum)\cdot(\trunckm\uuum -\uuum )\zeta.
\end{split}\end{align*}
This, \eqref{g4.1.14} and the product-rule leads to 
\begin{align}\begin{split}\label{g4.1.16}
&\int_{\omt} \Big(\grn(\nabla\om^M,\gre)-\widehat{\grn}\Big):(\nabla\trunckm\uuum)\zeta\\
&=\int_{\omt} (\widehat{\grn}-\grn(\nabla\om^M,\gre)):(\trunckm\uuum\otimes\nabla\zeta)\\
&\quad+\int_{\omt} \vek{H}_{2,M} :\nabla((\trunckm\uuum)\zeta)+\int_{\omt} \vek{k}_M\cdot (\trunckm\uuum)\zeta\\
&\quad+\frac{1}{2}\int_{\omt}(2\uuum\cdot \trunckm\uuum   -\abs{\trunckm\uuum}^2  )\partial_t \zeta\\
&\quad+\int_{E_{m,k}}(\partial_t \trunckm\uuum)\cdot(\uuum-\trunckm\uuum )\zeta\\
&=:\boldsymbol{1}_{M,k}+\boldsymbol{2}_{M,k}+\boldsymbol{3}_{M,k}+\boldsymbol{4}_{M,k}+\boldsymbol{5}_{M,k}.
\end{split}\end{align}
From now on we are able to use exactly the same ideas, which have been used in \cite{die-ru-wolf} to finish the proof. 
For the convenience of the reader we sketch them here. For any fixed $k \in \Nat$ we will pass to the 
limit in $M\rightarrow \infty$ in every integral of \eqref{g4.1.16} separately.
\begin{enumerate}
\item[(i)] $\limsup_{M\rightarrow \infty} (\abs{\boldsymbol{1}_{M,k}}+\abs{\boldsymbol{3}_{M,k}})=0$. \\
Since $\widehat{\grn}-\grn(\nabla\om^M,\gre)$ and $\vek{k}_M$ are bounded in $L^{p'}(\omt)$ by
 \eqref{g4.1.13}, we only need H\"{o}lder's inequality, the continuity result in Lemma 
 \ref{truncstetig} and \eqref{g4.1.13} to prove
\begin{align*}
\abs{\boldsymbol{1}_{m,k}}&\leq \norm{\zeta}_{L^{\infty}(\omt)} \norm{\widehat{\grn}-\grn(\nabla\om^M,\gre)}_{L^{p'}(\omt)}
 \norm{\trunckm\uuum}_{L^p(\omt)} \\
&\leq c \,\norm{\uuum}_{L^p(\omt)} \leq  c\, \norm{\uuum}_{L^{2\sigma}(\omt)} \rightarrow 0
\end{align*}
and
\begin{align*}
\abs{\boldsymbol{3}_{M,k}}&\leq \norm{\zeta}_{L^{\infty}(\omt)} \norm{\vek{k}_M}_{L^{p'}(\omt)}
 \norm{\trunckm\uuum}_{L^p(\omt)} \\
&\leq c\, \norm{\uuum}_{L^p(\omt)} \leq  c\, \norm{\uuum}_{L^{2\sigma}(\omt)} \rightarrow 0.
\end{align*}
\item[(ii)] $\limsup_{M\rightarrow \infty} \abs{\boldsymbol{2}_{M,k}}=0$. \\
We estimate
\begin{align*}\begin{split}
\abs{\boldsymbol{2}_{M,k}} \leq \norm{\vek{H}_{2,M}}_{L^1(\Omega_T)}
 \norm{\nabla((\trunckm\uuum)\zeta)}_{L^\infty(\Omega_T)}.
\end{split}\end{align*} 
The boundedness of $\Omega_T$ and \eqref{g4.1.13} implies $\vek{H}_{2,M}\rightarrow \bfzero$  in 
$L^1(\Omega_T)$, so that it remains to prove that for fixed $k \in \Nat$ the sequence
$(\nabla((\trunckm\uuum)\zeta))_{M\in\Nat}$ is bounded in $L^\infty(\operatorname{supp} (\zeta))$.
In view of 
\eqref{lmkabsch} 
we conclude that $\lMk$, $\aMk=\lMk^{2-p}$, $\lMk^{-1}$ and $\aMk^{-1}$ are all 
bounded in $M$ for fixed $k\in\Nat$. Additionally
\[ \inf_{M\in\Nat} \delta_{\aMk,K}=d_{\aMk}(K, \partial \Omega_T) >0 \]
holds for fixed $k \in \Nat$ as well. According to Theorem \ref{paralip}(ii) 
we are able to estimate
\begin{align*}\begin{split}
\norm{\nabla((\trunckm\uuum)\zeta)}_{L^\infty(K)} &\leq c(\lMk+\aMk^{-1} \delta_{\aMk,K}^{-3-3}\norm{\uuum}_{L^1(\emk)}   )\\
&\quad +c(1+\aMk^{-1} \delta_{\aMk,K}^{-3-2}\norm{\uuum}_{L^1(\emk)} )\\
&\leq c(k),
\end{split}\end{align*}
so that altogether $\limsup_{M\rightarrow \infty} \abs{\boldsymbol{2}_{M,k}}=0$ holds.
\item[(iii)] $\limsup_{M\rightarrow \infty} \abs{\boldsymbol{4}_{M,k}}=0$. \\
Using again H\"{o}lder's inequality, the continuity result in Lemma 
 \ref{truncstetig}  and \eqref{g4.1.13} we estimate
 \begin{align*}\begin{split}
\limsup_{M\rightarrow \infty} \abs{\boldsymbol{4}_{M,k}} \leq \limsup_{M\rightarrow \infty} c\,(1+ \norm{\partial_t \zeta}_{L^\infty(\Omega_T)})\norm{\uuum}_{L^2(\Omega_T)}^2 =0.
 \end{split}\end{align*}
\item[(iv)]$\limsup_{M\rightarrow \infty} \abs{\boldsymbol{5}_{M,k}}  \leq c \, 2^{-k}$.\\
Using $\operatorname{supp} \zeta=K$ and Theorem \ref{paralip} (iii) we estimate
\begin{align*}\begin{split}
\abs{\boldsymbol{5}_{M,k}} &\leq \norm{(\partial_t \trunckm\uuum)\cdot(\uuum-\trunckm\uuum )}_{L^1(\emk\cap K)}\\
&\leq c\, \aMk^{-1} \abs{\emk} (\lMk +\aMk^{-1} \delta_{\aMk,K}^{-3-3} \norm{\uuum}_{L^1(\emk)})^2.
\end{split}\end{align*}
Since $\aMk^{-1}$, $\abs{\emk}$, $\lMk$ and $\delta_{\aMk,K}^{-3-3}$ are all bounded in 
$M$ for fixed $k\in\Nat$ and since \[\norm{\uuum}_{L^1(\emk)} \leq \norm{\uuum}_{L^1(\Omega_T)} \stackrel{M\rightarrow \infty}{\longrightarrow}0\]
holds due to \eqref{emkTeilmenge} and \eqref{g4.1.13}, we conclude that
\begin{align*}\begin{split}
\limsup_{M\rightarrow \infty}\abs{\boldsymbol{5}_{M,k}} \leq \limsup_{M\rightarrow \infty} c\, \aMk^{-1} \abs{\emk} \lMk^2 = \limsup_{M\rightarrow \infty} c\,\lMk^p \abs{\emk} \stackrel{\eqref{emkAbsch}}{\leq}c \, 2^{-k}.
\end{split}\end{align*}
\end{enumerate}
Altogether, (i)--(iv) imply
\begin{align}\begin{split}\label{g5.1.1}
\limsup_{M\rightarrow \infty} \Big\vert\int_{\omt} \Big(\grn(\nabla\om^M,\gre)-\widehat{\grn}\Big):(\nabla\trunckm\uuum)\zeta\, \Big\vert \leq c\, 2^{-k}.
\end{split}\end{align}
On the other hand, from \eqref{SNbound}, \eqref{g4.1.13}, Theorem \ref{paralip} (ii)
and the boundedness of $\aMk^{-1}$, $\delta_{\aMk,K}^{-1}$ for fixed $k\in \Nat$, we conclude
\begin{align}\begin{split}\label{g5.1.2}
&\limsup_{M\rightarrow \infty} \Big\vert\int_{\emk} \Big(\grn(\nabla\om^M,\gre)-\widehat{\grn}\Big):(\nabla\trunckm\uuum)\zeta\, \Big\vert\\
&\quad\leq \limsup_{M\rightarrow \infty} \norm{\grn(\nabla\om^M,\gre)-\widehat{\grn}}_{L^{p'}(\omt)} \norm{\nabla\trunckm\uuum}_{L^\infty(K)}\abs{\emk}^{\frac{1}{p}}\\
&\quad\leq\limsup_{M\rightarrow \infty}  c\, \lMk \abs{\emk}^{\frac{1}{p}} \leq c\, 2^{-\frac{k}{p}}.
\end{split}\end{align}
Since $\nabla \trunckm\uuum =\nabla\uuum=\nabla\om^M-\nabla \om$ holds on $\omt\setminus\emk$ 
due to Definition 
\ref{truncdef}, the estimates \eqref{g5.1.1} and \eqref{g5.1.2} imply
\begin{align}\begin{split}\label{g5.1.3}
\limsup_{M\rightarrow \infty} \Big\vert\int_{\omt\setminus\emk} \Big(\grn(\nabla\om^M,\gre)-\widehat{\grn}\Big):(\nabla\om^M-\nabla \om)\zeta\, \Big\vert \leq c\, 2^{-\frac{k}{p}}.
\end{split}\end{align}
Due to \eqref{gmkabsch}, \eqref{fmkab}, \eqref{hmkabsch} and \eqref{g5.1.3} we are able to to find for each $k \in \Nat$ a number $M_k \in \Nat$ such that
\begin{align}\begin{split}\label{g5.1.4}
\Big\vert\int_{\omt\setminus E_{M_k,k}} &\Big(\grn(\nabla\om^{M_k},\gre)-\widehat{\grn}\Big):(\nabla\om^{M_k}-\nabla \om)\zeta\, \Big\vert \leq c\, 2^\frac{-k}{p},\\
\abs{G_{M_k,k}}&\leq c\,2^{-k},\quad
\abs{F_{M_k,k}}\leq c\,2^{-k},\quad
\abs{H_{M_k,k}}\leq c\,2^{-k}.
\end{split}\end{align} 
Now we set $\zeta_k:= \zeta \chi_{\omt\setminus E_{M_k,k}}$. Clearly we have 
\begin{align*}
\zeta_k \rightarrow \zeta \quad \text{pointwise in } \bigcup_{k=1}^\infty\bigcap_{\ell=k}^\infty
\big( \omt\setminus E_{M_\ell,\ell}\big)= \omt\setminus \bigcap_{k=1}^\infty 
\bigcup_{\ell=k}^\infty E_{M_\ell,\ell}.
\end{align*}
From \eqref{g5.1.4} we conclude that 
$\abs{\bigcap_{k=1}^\infty \bigcup_{\ell=k}^\infty E_{M_\ell,\ell}} =0$, so that 
\begin{align*}
\zeta_k\rightarrow \zeta \quad \text{holds a.e. in } \omt.
\end{align*}
This in turn implies the strong convergence of 
$\widehat{\grn}\zeta_k \rightarrow \widehat{\grn}\zeta$ in $L^{p'}(\omt)$ and 
$(\nabla\om)\zeta_k \rightarrow (\nabla\om)\zeta$ in $L^{p}(\omt)$. Now this together with \eqref{g5.1.4} implies
\begin{align*}
\lim\limits_{k\rightarrow \infty}  &\Big\vert\int_{\omt}  \grn(\nabla\om^{M_k},\gre): (\nabla\om^{M_k})\zeta_k - \widehat{\grn}: (\nabla \om)\zeta_k\,  \Big\vert\\
&\leq \lim\limits_{k\rightarrow \infty} \Big\vert\int_{\omt} \Big(\grn(\nabla\om^{M_k},\gre)-\widehat{\grn}\Big):(\nabla\om^{M_k}-\nabla \om)\zeta_k\, \Big\vert\\
&\quad+ \lim\limits_{k\rightarrow \infty} \Big\vert\int_{\omt} \grn(\nabla\om^{M_k},\gre): (\nabla\om)\zeta_k         
+ \widehat{\grn}: (\nabla \om^{M_k,k})\zeta_k - 2\widehat{\grn}:(\nabla\om)\zeta_k       \, \Big\vert\\
&=0,
\end{align*}
so that we have
\begin{align*}
\lim\limits_{k\rightarrow \infty} \int_{\omt}  \grn(\nabla\om^{M_k},\gre): (\nabla\om^{M_k})\zeta_k
= \int_{\omt} \widehat{\grn}: (\nabla \om)\zeta.
\end{align*}
The local Minty trick (cf.~\cite[Lemma A.2]{Wolf}) implies 
$\widehat{\grn}\zeta= \grn(\nabla\om,\gre)\zeta$ a.e.~in $\omt$. 
\end{Bew}
\begin{Bem}
It is also possible to derive the weak continuity of the solutions $\vek{v}$ and $\om$. 
Due to $\eqref{konvmain}$ and $\eqref{einbettung1}$ we have that 
$\vek{v} \in C(0,T;V_{\sigma'}(\Omega)^\ast)$ and $\om \in C(0,T;W_0^{1,\sigma'}(\Omega)^{\ast})$. 
Since the solutions belong to the spaces $L^\infty(0,T;H(\Omega))$ and $L^\infty(0,T;L^2(\Omega))$,
respectively, we are able to conclude (cf.~\cite[Lemma II.5.9]{boy-fab}) that $\vek{v}$ is weakly continuous with 
values in $H(\Omega)$, while $\om$ is weakly continuous with values in $L^2(\Omega)$. Let us have 
a quick look why the equality $\vek{v}(0)=\vek{v}_0$ also holds in the sense of this weak continuity. 
Since the approximative solutions $\vek{v}^M$ posses the required initial values, we are able to 
derive from the integration by parts formula in $\eqref{partielleInteg}$ that
\begin{align*}\begin{split}
(\vek{v}_0 ,\fett{\eta})_H=\int_0^T \langle \frac{d\vv^M}{dt},\fett{\eta} \rangle_{V_{\sigma'}}\zeta
+ \int_{\omt} \vv^M\cdot \fett{\eta}\partial_t\zeta
\end{split}\end{align*}
holds for any $\fett{\eta} \in \mathcal{V}(\Omega)$ and $\zeta \in C^\infty([0,T])$ with 
$\zeta(0)=-1$ and $\zeta(T)=0$. With this $\zeta$ we have 
$\zeta\vek{v} \in W^{1,p,\sigma}(0,T;V_p(\Omega), V_{\sigma'}(\Omega)^{\ast})$ which implies
\begin{align*}\begin{split}
(\vek{v}(0) ,\fett{\eta})_H=\langle\vek{v}(0) ,\fett{\eta}\rangle_{V_{\sigma'}}= \int_0^T \langle \frac{d\vv}{dt}
,\fett{\eta} \rangle_{V_{\sigma'}}\zeta+ \int_{\omt} \vv\cdot \fett{\eta}\partial_t\zeta,
\end{split}\end{align*} 
so that we get $(\vek{v}_0 ,\fett{\eta})_H=(\vek{v}(0) ,\fett{\eta})_H$ due to \eqref{konvmain}. 
This implies $\vek{v}_0=\vek{v}(0)$ in $H(\Omega)$. The argumentation for $\om_0=\om(0)$ in 
$L^2(\Omega)$ is the same.

\end{Bem}

\section*{References}
\def\cprime{$'$} \def\cprime{$'$} \def\cprime{$'$}
\ifx\undefined\bysame
\newcommand{\bysame}{\leavevmode\hbox to3em{\hrulefill}\,}
\fi


\end{document}